\documentclass[a4paper]{article}
\usepackage{amsfonts,amsmath,amssymb}
\usepackage{graphicx}
\usepackage{bm}
\input{epsf.sty}
\usepackage{subfigure}
\usepackage{epsfig}
\usepackage{footmisc}
\usepackage{url}
\usepackage{amsthm}
\usepackage{amscd,enumerate,latexsym}
\usepackage{multirow}
\usepackage{cite}
\usepackage{color}
\usepackage{hyperref}
\usepackage{float}
\usepackage{algorithm,algorithmic}
\usepackage{indentfirst}
\usepackage{threeparttable}
\usepackage{booktabs}
\usepackage{rotating}
\usepackage{placeins}
\usepackage{breakurl}
\newtheorem{theorem}{\bf Theorem}[section]

\newtheorem{proposition}{Proposition}[section]

\newcommand{\ms}[1]{\medskip}
\newcommand{\bs}[1]{\bigskip}
\newcommand{\ssk}[1]{\smallskip}

\begin{document}
\begin{center}
\Large\bf{A Hessenberg-type Algorithm for Computing PageRank Problems}\\

\quad\\
\normalsize Xian-Ming Gu\footnote[1]{School of Economic Mathematics/Institute of
Mathematics, Southwestern University of Finance and Economics, Chengdu, Sichuan
611130, P.R. China and Bernoulli Institute for Mathematics, Computer Science and Artificial Intelligence,
University of Groningen, Nijenborgh 9, P.O. Box 407, 9700 AK Groningen, The Netherlands. E-mail: {\tt guxianming@live.cn,~guxm@swufe.edu.cn}.},~
Siu-Long Lei\footnote[2,*]{Corresponding author: Department of Mathematics, Faculty of Science and Technology,
University of Macau, Avenida da Universidade, Macao, P.R. China. E-mail: {\tt sllei@um.edu.mo}.},~
Ke Zhang\footnote[3]{College of Arts and Sciences, Shanghai Maritime University,
Shanghai 201306, P.R. China. E-mail: {\tt xznuzk123@126.com}.},~
Zhao-Li Shen\footnote[4]{Department of Applied Mathematics, College of Science, Sichuan Agricultural
University, Yaan, Sichuan 625014, P.R. China. E-mail: {\tt szlxiaoyao@163.com}.},\\
Chun Wen\footnote[5]{School of Mathematical Sciences, University of Electronic
Science and Technology of China, Chengdu, Sichuan 611731, P.R. China. E-mail: {\tt
wchun17@163.com}.},~
Bruno Carpentieri\footnote[6]{Facolt\`{a} di Scienze e Tecnologie informatiche,
Libera Universit\`{a} di Bolzano, Dominikanerplatz 3 - piazza Domenicani, 3 Italy
- 39100, Bozen-Bolzano. E-mail: {\tt Bruno.Carpentieri@unibz.it}.}
\end{center}

\vspace{0.4cm}

\mbox{\bf Abstract:}~~PageRank is a widespread model for analysing the relative relevance of nodes within large graphs arising in several applications.
In the current paper, we present a cost-effective Hessenberg-type method built
upon the Hessenberg process for the solution of difficult PageRank
problems. The new method is very competitive with other popular algorithms in this field, such as Arnoldi-type methods, especially when the damping factor is close to $1$ and the dimension of the search subspace is large. The convergence and the complexity of the proposed algorithm are investigated. Numerical experiments are reported to show the efficiency of the new solver for practical PageRank computations.

\mbox{\bf Keywords:} PageRank vector, Hessenberg process, Arnoldi process, Krylov subspace
method, Ritz values.

\mbox{\bf AMS classifications: } 65F15, 65F10, 65Y20.

\section{Introduction}
\setcounter{equation}{0}
The PageRank model was originally introduced by S.~Brin and L.~Page in 1999~\cite{Page1998} to develop fast web search engines, and then studied and enhanced in a vast number of research papers (see, e.g.,~\cite{Kamvar03,Kamvar04x,Langv05,Langville05,Berkhin05,LanMey2006}). The model provides a powerful network centrality measure to identify the most important nodes within large graphs arising in several applications fields, such as in chemistry, bioinformatics, neuroscience, bibliometrics, and others~\cite{Gleich15}. In the original Web problem, the PageRank algorithm determines the ranking of each Web page by computing the stationary probability vector of a random walking process on the Web link graph, which is a directed graph representing the linking structure of the Web~\cite{LanMey2006,Bryan06}. The Web link graph is a binary matrix $G\in\mathbb{N}^{n\times n}$, where $n$ denotes the number of pages, such that $G(i,j)=1$ when page $j$ has a link pointing to page $i$, and $G(i,j)=0$ otherwise. From a linear algebra viewpoint, the algorithm finds the vector $\bm x$ that satisfies
\begin{equation}\label{eq:pagerank}
	A{\bm x} = {\bm x}, \quad \|{\bm x}\|_1 = 1,~~{\bm x} > 0,
\end{equation}
that is it finds the principal unit positive eigenvector $\bm{x}$~\cite{Bryan06} of the Google matrix
\begin{equation}
A = \alpha (P + {\bm v}{\bm d}^{\top}) + (1 - \alpha){\bm v}{\bm e}^{\top}.
\label{eq1.1}
\end{equation}
In Eq.~(\ref{eq1.1}), matrix $P\in \mathbb{R}^{n\times n}$ is called the \textit{transition matrix} with respect to the random walking process, and is defined as
\begin{equation}
P(i,j)=
\begin{cases}
  \frac{1}{\sum\limits^{n}_{k=1} G(k,j)}, & \mbox{if } G(i,j)=1, \\
  0, & \mbox{otherwise}.
\end{cases}
\label{P}
\end{equation}
The \textit{damping factor} $0<\alpha<1$ defines the probability that random
Web surfers choose a random link from the page they are visiting~\cite{Wu2013x}.
The \textit{teleporting vector} ${\bm v}=[v_1,\cdots,v_n]^\top\in \mathbb{R}^{n\times 1}$ (${\bm v}\geq 0$ and $\|{\bm v}\|_1=1$) defines the probability $v_i$ that the Web surfer jumps to an external page $i$. Finally, ${\bm d}\in \mathbb{N}^{n\times 1}$ is such that ${\bm d}(i)=1$ if page $i$ has no hyperlink and $0$ otherwise, and ${\bm e}=[1,1,\ldots,1]^{\top}\in \mathbb{N}^{n\times 1}$.

The value of the damping factor $\alpha$ plays an important role in the PageRank model. Theoretically, it represents an upper bound $0 < |\lambda_2| \leq \alpha < 1$ for the second largest eigenvalue, $\lambda_2$, of $A$. Further properties of the Google matrix can be found in~\cite{Kamvar03,Langv05,Langville05,Berkhin05,LanMey2006,Cicone10}. For low values of $\alpha$ (e.g., $\alpha=0.85$), $\lambda_2$ is well separated from the largest eigenvalue of $A$, which is $\lambda_1=1$, ensuring rapid convergence of the power algorithm applied to problem~(\ref{eq:pagerank}). On the other hand, convergence tends to slow down noticeably when $\alpha$ is very close to $1$, requiring more robust algorithms than the simple power method. Some computational approaches proposed in the literature include Monte Carlo methods~\cite{Avrach07}, adaptive algorithms~\cite{Kamvar04x,Liu2015}, extrapolation techniques~\cite{Kamvar03,LanMey2006,Tan2017}, singular value decompositions~\cite{Jia1997,Jia1999,Jia2000}
reordering~\cite{Langville06,Lin2009} and inner-outer solution strategies~\cite{Gleich10}.

A significant amount of work has been devoted in the last years to the use of Krylov subspace methods based on the Arnoldi decomposition~\cite{Heyouni98,Saad2010} for large PageRank computations, mainly due to their memory efficiency and attractive inherent parallelism. Golub and Greif extended the refined Arnoldi procedure to PageRank by forcing a relevant shift to be 1, being able to circumvent the drawbacks due to complex arithmetic and showing overall very good algorithmic efficiency~\cite{Golub06}. Many techniques have attempted to combine the conventional Arnoldi method and the power algorithm to produce faster solvers, e.g., the Power-Arnoldi~\cite{Wu2007,Yin2010,Wu2013x},
the Arnoldi-Extrapolation~\cite{Wu2010} and the Arnoldi-Inout \cite{Gu2017a} methods. In the technique proposed in~\cite{Yin2012}, the weighted least squares problem
is changed adaptively according to the component of the residual. Then, the generalized Arnoldi method is used to compute the approximate PageRank vector. However, when the dimension of the Krylov subspace is large, Arnoldi-based solvers tend to become very expensive in terms of memory and computational costs; on the other hand, if the dimension of the Krylov subspace is low, they sometimes fail to accelerate the basic power method, especially when the damping factor is high~\cite{Golub06,Wu2007,Yin2010,Wu2013x}. Similarly to the restarted GMRES algorithm~\cite{WuWanJin2012}, they may stagnate in many circumstances~\cite{Wu2010y}.

Motivated by costs considerations, other work developed PageRank solvers based on the Bi-Lanczos orthogonalization procedure~\cite[pp. 139-145]{Saad2010} (see e.g.~\cite{Freund94,Teramoto}) instead of Arnoldi. In this paper, we look in particular at the Hessenberg reduction
process~\cite{Hessen40,Wilkinson,Sadok1999,Stephens99} that was introduced by K.~Hessenberg in 1940~\cite{Hessen40}, and revived recently to establish a number of cost-effective Krylov subspace solvers for sparse matrix systems, due to its lower arithmetic and storage requirements. The method has been extended to compute the characteristic polynomial of matrices~\cite{Hessen40,Wilkinson,House2010}, to solve general nonsymmetric linear systems~\cite{Heyouni98,Sadok1999,Stephens99,Heyouni08,Sadok2012}, including systems with multiple right-hand sides
\cite{Heyouni01,Heyouni05x,Zhang14x,Amini18x,Amini18y}
and multi-shifted coefficient matrices~\cite{Gu2018JCAM,Gu2018X,Ramezani18}, other types of matrix equations~\cite{Heyouni05x,Addam2017,Heyouni19},
the action of a matrix function $f(A){\bm v}$ \cite{Ramezani18} and other related problems~\cite{Heyouni99x}. Theoretical and numerical studies have investigated the mathematical properties of the Hessenberg process, especially in relation with the more conventional Arnoldi procedure. The Arnoldi method was first introduced in 1951 as a means
of reducing a dense matrix $A$ into a Hessenberg form by unitary transformations, whereas the Hessenberg
process applies \textit{similarity transformations}~\cite{Businger} and is more suitable for parallel computing. In his paper, Arnoldi hinted that the eigenvalues of the Hessenberg matrix obtained after $k\ll n$ steps, where $n$ is the size of $A$, could provide accurate approximations of some eigenvalue of $A$. It was later discovered that this
strategy can lead to efficient techniques for approximating eigenvalues of large sparse matrices. In
the current work, following a similar development, we  modify the Hessenberg process to establish a new family of eigenvalue solvers. We combine the new solvers with the refined and explicitly restarted techniques introduced in \cite{Jia1997,Golub06} to compute realistic PageRank problems. Finally, we analyse their convergence behavior and computational complexity.

The rest of this paper is organized as follows. In Section \ref{sec2}, the
Hessenberg process is introduced and a novel family of eigenvalue solvers
based on this procedure is described. Moreover, theoretical aspects of such
eigenvalue solvers are highlighted in comparison with the classical
Arnoldi-like methods. In Section~\ref{sec3}, we derive the Hessenberg-type method
with explicit restarting and refined techniques for computing PageRank. Both the
convergence behavior and the computational cost of the proposed method are discussed.
Numerical results in Section \ref{sec4} show the effectiveness of the
proposed algorithm, also against other popular PageRank algorithms. In Section
\ref{sec5}, we present some conclusions arising from our study.

\section{The Hessenberg process with applications to eigenvalue computations}
\label{sec2}
In this section, we briefly review the Hessenberg procedure that is at the basis
of our development. We recall some fundamental properties of the algorithm and then we describe a Hessenberg-based projection technique for computing eigenvalues of large nonsymmetric matrices. Our theoretical analyses demonstrates the feasibility of the method, showing some computational advantages over the more conventional Arnoldi procedure.
\subsection{The Hessenberg process}
The Hessenberg process is an an oblique projection technique that reduces
a given nonsymmetric matrix $A\in\mathbb{R}^{n\times n}$
to a Hessenberg form~\cite[pp. 377-381]{Wilkinson,Gu2018JCAM}. Originally, the method was described as a way to compute the characteristic polynomial of a matrix~\cite{Hessen40}. The basic procedure is presented  in Algorithm~\ref{Ad-SGMRES}, where a pivoting strategy is included to ensure numerical stability.

\begin{algorithm}[!htbp]
\caption{The Hessenberg procedure with pivoting strategy}
\begin{algorithmic}[1]
\REQUIRE a general (nonsymmetric) square matrix $A\in\mathbb{R}^{n\times n}$, an initial vector
${\bm q}_0$ and an integer $m$
\ENSURE $L_k=\left[ {\bm l}_1,\ldots,{\bm l}_m \right]\in
\mathbb{R}^{n\times m}$ and $H_m\in\mathbb{R}^{(m+1)\times m}$
\STATE Initialize the permutation vector ${\bm p} = [1,2,\ldots,n]^\top$ and determine $i_0$ such that $|({\bm v})_{i_0}|= \|{\bm v}\|_{\infty}$
\STATE Compute $\beta = ({\bm v})_{i_0}$, then ${\bm l}_1 = {\bm v}/\beta$ and ${\bm p}(1)
\leftrightarrow {\bm p}(i_0)$, where ``$\leftrightarrow$" is used to swap contents \cite{Sadok1999}.
\FOR{$j = 1,2,\ldots,m$,}
\STATE Compute ${\bm u} = A{\bm l}_j$
\FOR{$i = 1,2,\ldots,j$,}
\STATE $h_{i,j} = ({\bm u})_{{\bm p}(i)}$
\STATE ${\bm u} = {\bm u} - h_{i,j}{\bm l}_i$
\ENDFOR
\IF{($j < n$ and ${\bm u} \neq {\bm 0}$)}
\STATE Determine $i_0 \in \{j + 1,...,n\}$ such that $|({\bm u})_{{\bm p}(i_0)}|=\|({\bm u})_{{\bm p}(j+1):{\bm p}(n)}\|_{\infty}$;
\STATE $h_{j + 1,j} = ({\bm u})_{{\bm p}(i_0)}$;~~${\bm l}_{j+1} = {\bm u}/h_{j + 1,j}$;~~${\bm p}(j+1)\leftrightarrow {\bm p}(i_0)$
\ELSE
\STATE $h_{j + 1,j} = 0$; \textbf{Stop} \hfill {\tt $\rhd$ Happy breakdown}
\ENDIF
\ENDFOR
\end{algorithmic}
\label{Ad-SGMRES}
\end{algorithm}
Let $L_k=\left[ {\bm l}_1,\ldots,{\bm l}_m \right]$ denote a matrix,
$\bar{H}_m=\left[ h_{i,j} \right]$ be an upper Hessenberg matrix and $H_m$ the submatrix obtained from $\bar{H}_m$ by deleting its last row. Finally, we denote by $\mathcal{P}^{{\color{red}\top}}_k = [{\bm e}_{p_1},{\bm e}_{p_2},\ldots,{\bm e}_{p_n}]$ where the scalars $p_i$'s (for $i = 1,\ldots,n$) are defined
in Algorithm~\ref{Ad-SGMRES}. After $k$ steps of Algorithm~\ref{Ad-SGMRES}, the following matrix equation can be easily established,
\begin{equation}
\begin{split}
AL_k & = L_{k+1}\bar{H}_k \\
     & = L_kH_k + h_{k+1,k}{\bm l}_{k+1} {\bm e}^{\top}_k,
\end{split}
\label{eq1.7}
\end{equation}
and $\mathcal{P}_kL_k$ is lower trapezoidal~\cite{Sadok1999,Heyouni08}. Unlike Arnoldi, the Hessenberg procedure with pivoting  is not guaranteed to be backward stable in finite precision arithmetic~\cite{Businger}. However, the backward error is reported to be small for most practical problems~\cite{Stephens99,Heyouni08,Gu2018JCAM}. This is also confirmed by our computational experiences.
We did not observe noticeable numerical instabilities due to the non-orthogonality of the Krylov basis in our numerical experiments; see also~\cite{Gu2018JCAM} for a discussion of this topic.
\subsection{Approximation of eigenpairs based on the Hessenberg process}
\label{sec2.2}
Methods to approximate eigenpairs of a large nonsymmetric matrix $A$ usually compute them from the Hessenberg decomposition of $A$ given by Eq.~(\ref{eq1.7}).
The upper Hessenberg matrix $H_m$ can be seen as the projection of $A$ onto the Krylov subspace
\begin{equation}
\mathcal{K}_m(A,{\bm v}) = {\rm span}\{{\bm v},A{\bm v},\ldots,A^{m-1}{\bm v}\},
\end{equation}
and the columns of the matrix $L_m$ are a basis of $\mathcal{K}_m(A,{\bm v})$.
Under certain conditions, the eigenvalues of $H_m$  converge to the eigenvalues of $A$~\cite{Wilkinson,Businger,Saad2010}. Various eigenvalue solvers are built upon this idea, differing each other mainly on the type of projection technique that is used to derive the decomposition~(\ref{eq1.7}), for example the Arnoldi process, the Bi-Lanczos procedure and the Induced Dimension Reduction (IDR) strategy~\cite{Gutknecht,Astudillo}. The approximate eigenpairs of $A$ are retrieved in the form $(\theta_i,~{\bm x}^{(i)} = L_m{\bm y}^{(i)})$, where $\left(\theta_i, {\bm y}^{(i)}\right)$ are eigenpairs of the small dimensional matrix $H_m$, such that
\begin{equation}
H_m{\bm y}^{(i)} = \theta_i{\bm y}^{(i)}\quad {\rm with}~~\|{\bm y}^{(i)}\| = 1,~~i=1,2,\cdots,m.
\label{xxxku}
\end{equation}
A bound on the residual error can be established directly from Eq.~(\ref{eq1.7}), by writing
\begin{equation*}
\begin{split}
A{\bm x}^{(i)} - \theta_i{\bm x}^{(i)} & = AL_m{\bm y}^{(i)} - \theta_i L_m{\bm y}^{(i)}\\
& = h_{m + 1,m}{\bm l}_{m+1}{\bm e}^{\top}_m{\bm y}^{(i)}.
\end{split}
\end{equation*}
If we denote as $[{\bm y}^{(i)}]_m$ the $m$th component of the vector ${\bm y}^{(i)}$, then we obtain
\begin{equation}
\|A{\bm x}^{(i)} - \theta_i{\bm x}^{(i)}\| \leq |h_{m + 1,m}|\|{\bm l}_{m + 1,m}\|
\left|[{\bm y}^{(i)}]_m\right|
\end{equation}
or, if we normalize the vector ${\bm l}_m$,
\begin{equation}
\|A{\bm x}^{(i)} - \theta_i{\bm x}^{(i)}\| \leq |h_{m + 1,m}|
\left|[{\bm y}^{(i)}]_m\right|.
\label{eq2.6}
\end{equation}
This analysis is in line with the results described in~\cite{Jia1997}.

\begin{figure}[t]
\centering
\includegraphics[width=0.49\textwidth]{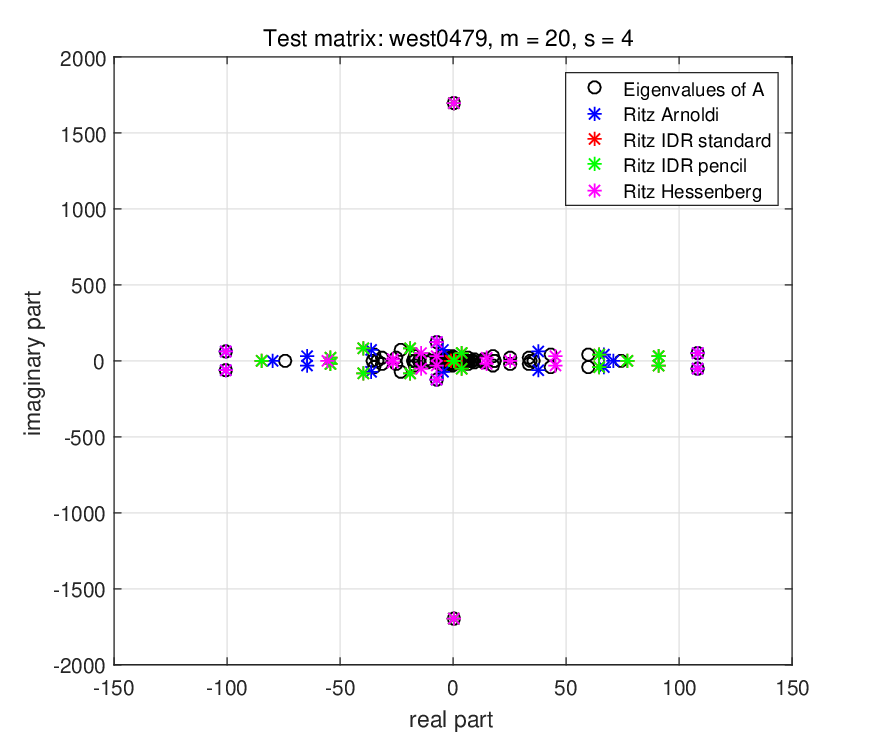}
\includegraphics[width=0.49\textwidth]{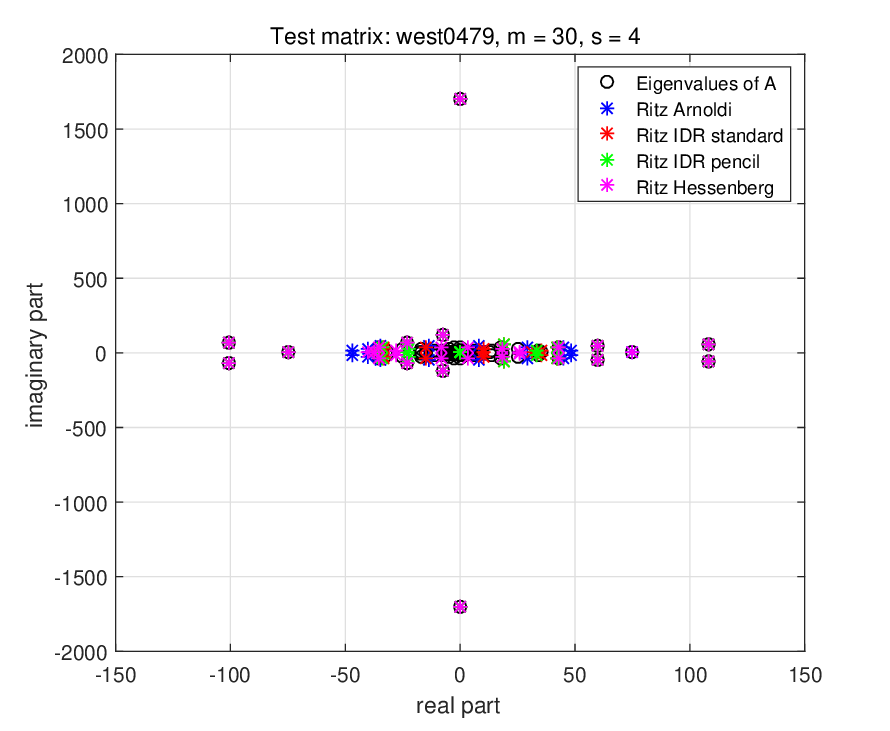}
\includegraphics[width=0.49\textwidth]{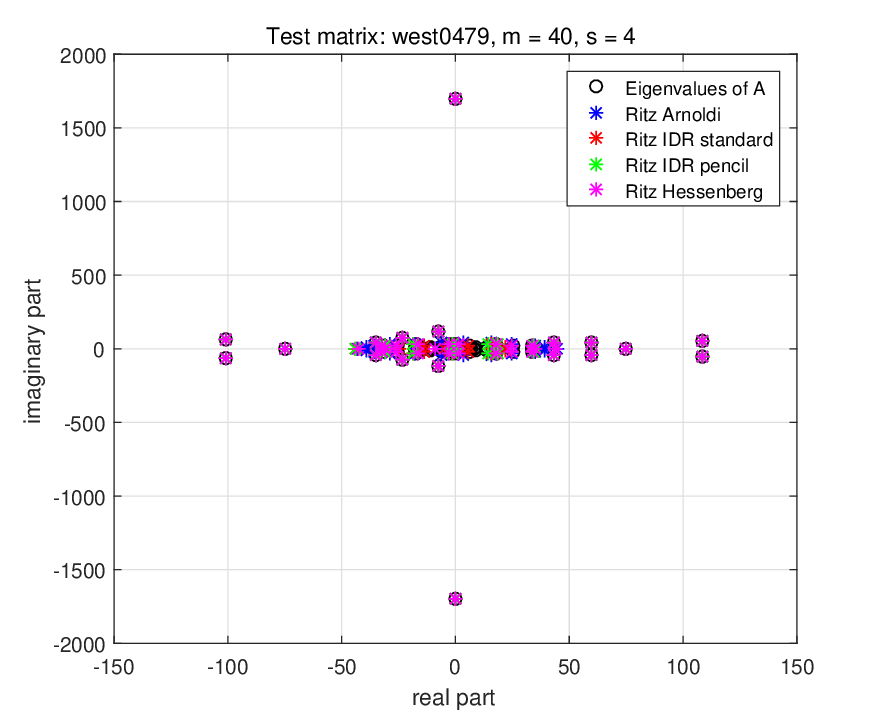}
\includegraphics[width=0.49\textwidth]{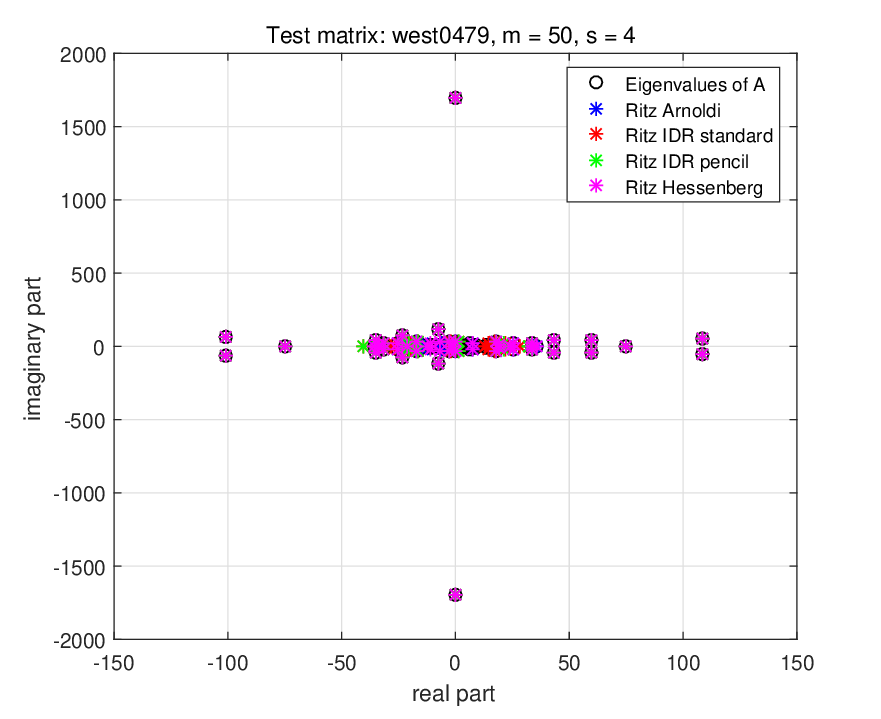}
\caption{Plots of the Ritz values generated by the IDR($s = 4$) factorization, the Sonneveld pencil, the Arnoldi and the Hessenberg procedures.}
\label{fig1}
\end{figure}

\begin{figure}[!htpb]
\centering
\includegraphics[width=0.49\textwidth]{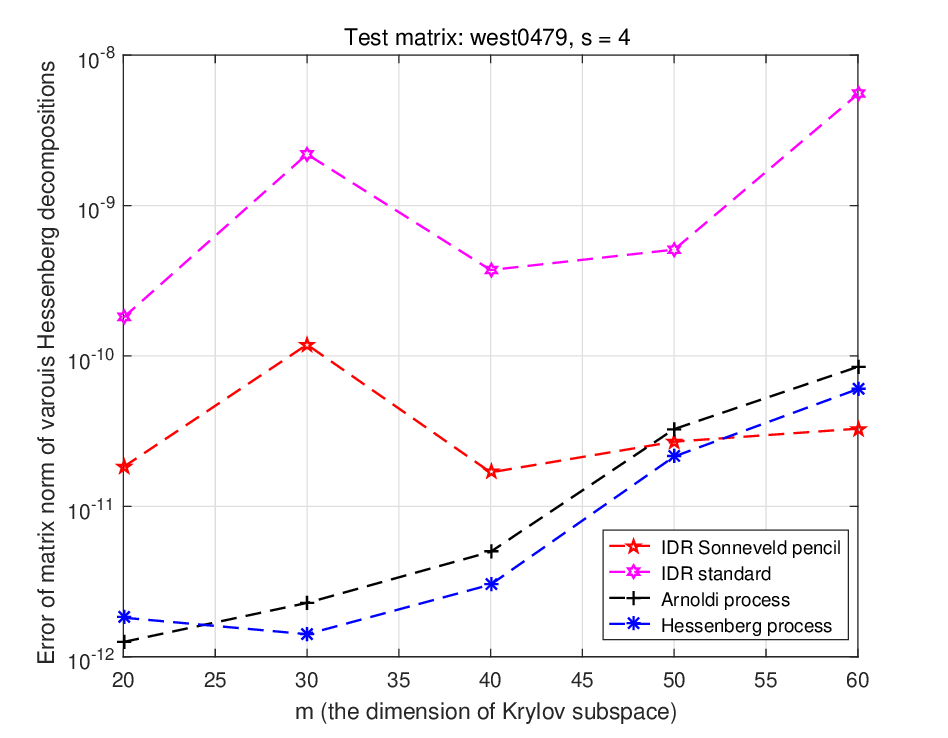}
\includegraphics[width=0.49\textwidth]{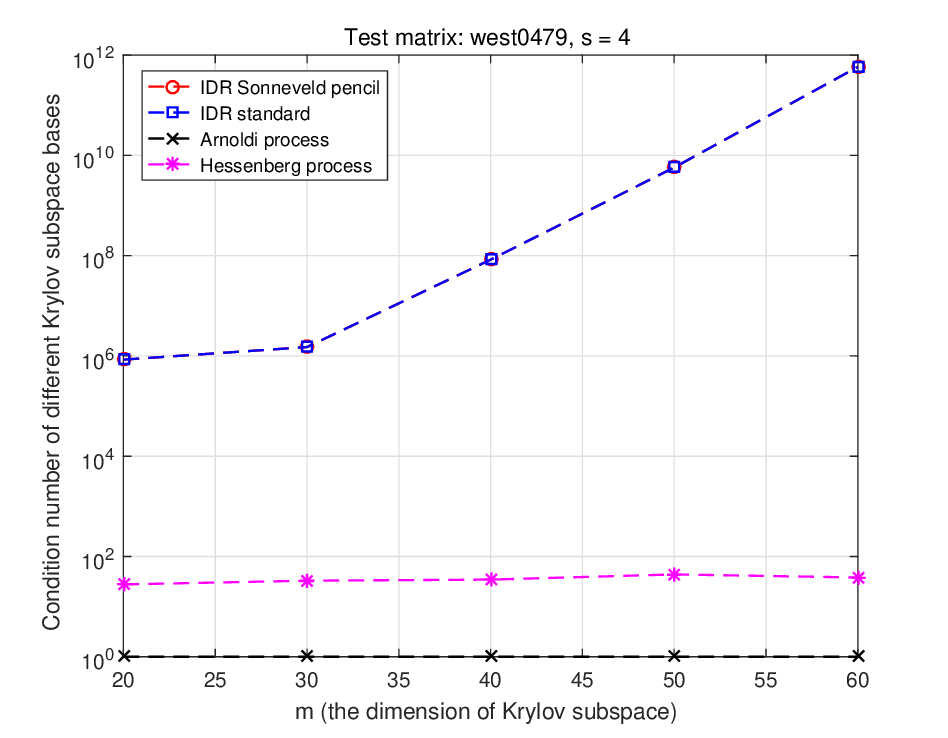}
\caption{The quality of different basis matrices; Left: the error of different Hessenberg decompositions;
Right: the condition number of basis matrices generated by different Hessenberg decompositions.}
\label{fig2}
\end{figure}

For the ``\texttt{west0479}" problem, a real-valued 479-by-479 sparse matrix that has both real and complex eigenvalues, in Fig.~\ref{fig1} we plot the eigenvalues of $A$ computed by the MATLAB command \texttt{eig}, and those of the Hessenberg matrix produced by the IDR($s = 4$) projection technique~\cite{Astudillo}, the Sonneveld pencil~\cite{Gutknecht,Astudillo}, the Arnoldi and the Hessenberg procedures. We clearly see that the Hessenberg procedure can estimate the exterior Ritz values very well, in some cases even slightly more accurately than the Arnoldi procedure. The condition number of the Krylov basis matrix $L_m$ is an effective metric to determine the accuracy of the method used. Fig.~\ref{fig2} illustrates that for the Hessenberg process, this condition number does not grow significantly when the dimension of the Krylov subspace increases. The numerical error of the Hessenberg decomposition often remains small~\footnote{Here we give its executable codes in the website: \url{https://github.com/Hsien-Ming-Ku/PageRank-Hessenberg}.}. This observation is also exemplified by the stochastic analysis presented in~\cite[Section 3.3]{Astudillo},
that can produce (numerical) evidence that the Hessenberg process can be as efficient as the Arnoldi process for practical eigenvalue computations.
\subsection{Relation between the Arnoldi and Hessenberg decompositions}
In this subsection, we provide more theoretical background supporting the use of the Hessenberg process to approximate effectively eigenpairs of a given nonsymmetric matrix $A$. The starting point of our analysis is a comparison between the Hessenberg decompositions computed by the Arnoldi and by the Hessenberg procedures. After $m$ steps of the Arnoldi method applied to $A$, starting with an initial vector ${\bm v}_0$ and assuming no breakdown, the following Hessenberg decomposition is derived:
\begin{equation}
	AV_m = V_mH_m + h_{m+1,m}{\bm v}_{m + 1}{\bm e}^{\top}_m = V_{m + 1}\bar{H}_m.
	\label{eq2.7}
\end{equation}

On the other hand, after $m$ steps of the Hessenberg procedure applied to $A$ with same initial vector ${\bm l}_0$, the resulting matrix factorization writes as
\begin{equation}
AL_m = L_mH^{(h)}_m + h^{(h)}_{m+1,m}{\bm l}_{m+1}{\bm e}^{\top}_m = L_{m + 1}\hat{H}^{(h)}_m.
\label{eq2.8}
\end{equation}
Differently from Arnoldi, however, the columns of $L_m$ in Eq.~\ref{eq2.8} are not mutually orthogonal. By computing the reduced QR factorization of $L_{m + 1}$,
\begin{equation}
L_{m + 1} = Q_{m+1}R_{m + 1},
\label{FAdSGMRESDReig}
\end{equation}
we can establish the following relation between Eq.~(\ref{eq2.8}) and Eq.~(\ref{FAdSGMRESDReig}):
\begin{equation}
AQ_m = Q_{m + 1}R_{m+1}\hat{H}^{(h)}_mR^{-1}_m.
\label{eq2.10}
\end{equation}
Due to the uniqueness of the Arnoldi decomposition (see Section 3.3 in \cite{Astudillo}), by comparing Eq.~(\ref{eq2.7}) and Eq.~(\ref{eq2.10}) we conclude that $Q_{m + 1} = V_{m+1}$ and
\begin{equation}
\bar{H}_m = R_{m+1}\hat{H}^{(h)}_m R^{-1}_m.
\end{equation}
The above result is summarized in the following proposition:
\begin{proposition} It follows that
\begin{eqnarray}
H_m & = & R_mH^{(h)}_mR^{-1}_m + \frac{{h^{(h)}_{m+1,m}}}{r_{m,m}}\tilde{{\bm r}}
{\bm e}^{\top}_m,\label{eq2.12}\\
\frac{h^{(h)}_{m+1,m}}{r_{m,m}} & = & \frac{h_{m+1,m}}{r_{m + 1,m + 1}}, \label{eq2.13}
\end{eqnarray}
where $\tilde{{\bm r}} = [r_{i,m+1}]^{m}_{i=1}$ is the vector containing the first
$m$ components of the $(m + 1)$th column of $R_{m + 1}$.
\label{pro2.1}
\end{proposition}
\noindent In fact, Eqs.~(\ref{eq2.12})--(\ref{eq2.13}) can be also found in~\cite{Sadok2012,Astudillo}. According to Eq.~(\ref{eq2.13}), in exact arithmetic both procedures produce upper Hessenberg matrices with the same eigenvalues. If the Arnoldi process terminates successfully (i.e., a happy breakdown situation $h_{m + 1,m}= 0$ occurs), so does the Hessenberg procedure ($h^{(h)}_{m + 1,m}= 0$). On the other hand, a direct consequence of Proposition~\ref{pro2.1} is that the Ritz values produced by the Arnoldi and by the Hessenberg processes are not the same. The condition number of the Krylov basis matrix $L_m$, which is the same as the condition number of matrix $R_m$, gives a clear indication of the accuracy of the eigenvalues of $H^{(h)}_m$ compared to those resulting from the Arnoldi process.

In conclusion, it can be expected that the Hessenberg process can produce feasible approximations of eigenpairs of large nonsymmetric matrices matrix.

\section{A Hessenberg-type algorithm for computing PageRank}\label{sec3}
In this section, we will propose a Hessenberg-based algorithm to compute the PageRank
vector, that is the positive unit eigenvector corresponding to the largest eigenvalue
of the Google matrix. Golub and Greif suggested that the explicitly restarted Arnoldi
process for computing eigenvalues and  eigenvectors should be implemented in complex arithmetic, thus it is not suitable (it needs to be  refined) to compute the PageRank vector~\cite{Golub06}.
Moreover, since the Hessenberg
process is similar in nature to Arnoldi, except that it produces a non-orthogonal basis of the Krylov subspace, we follow a similar development to the {\it refined Arnoldi method} for PageRank problems. In other words, the solver described in this section may be called the \textit{refined Hessenberg method} for PageRank. However, we do not approximate the eigenvectors of $A$ from those of $H_{m}$ (the so-called Ritz-like vectors). Instead, we compute the refined
Ritz-like vectors, i.e., the singular vectors associated with the smallest singular values of $A - \theta_iI$~\cite{Jia1997,Golub06}, where $\{\theta_i\}^{m}_{i=1}$ are named the Ritz-like values; cf. Eq. (\ref{xxxku}). The Hessenberg-based method enjoys similar numerical properties to the Arnoldi-based variant: an effective separation of the eigenvectors is ensured, complex arithmetic is avoided by using a shift equal to $1$ (due to the fact that the largest eigenvalue of the Google matrix is known), the smallest singular value converges more smoothly to zero than the largest Ritz value to $1$ \cite[Section 3]{Golub06}. The Hessenberg-type method for computing the PageRank vector is presented in Algorithm~\ref{Alg2}. The following convergence result can be established after each cycle of $m$ iterations of Algorithm~\ref{Alg2}.

\begin{algorithm}[!htbp]
\caption{The Hessenberg-type algorithm for computing PageRank}
\begin{algorithmic}[1]
\REQUIRE the Google matrix $A\in\mathbb{R}^{n\times n}$, initial guess
${\bm q}_0$ and $m$
\ENSURE the approximate PageRank vector ${\bm q}$
\STATE Choose $q_0$ and determine $i_0$ such that $|({\bm q_0})_{i_0}|= \|{\bm q_0}\|_{\infty}$, $m$
is the dimension of search subspace
\STATE Compute $\beta = ({\bm q_0})_{i_0}$, then ${\bm q}_0 = {\bm q_0}/\beta$
\FOR{$\ell = 1,2,\ldots,$ until convergence}
\STATE Compute $L_m$ and $H_{m+1,m}$ via running Algorithm \ref{Ad-SGMRES} with the inputs
$A,~m,~{\bm q}_0$
\STATE Compute the singular value decomposition: $H_{m+1,m} - [I_m;{\bm 0}] = U\Gamma V^\top$
\STATE Compute ${\bm q} = L_m{\bm v}_m$
\STATE Compute ${\bm r} = \sigma_mL_{m+1}{\bm u}_m$
\IF{$\|{\bm r}\|_1/\|{\bm q}\|_1 < tol$}
\STATE {\em Stop and exit}
\ELSE
\STATE Set ${\bm q}_0 = {\bm q}$
\ENDIF
\ENDFOR
\end{algorithmic}
\label{Alg2}
\end{algorithm}

\begin{theorem}
Let $Q_m = [{\bm q}_1,{\bm q}_2,\cdots,{\bm q}_m]$ be the matrix obtained from running $m$-steps of either the Arnoldi or the Hessenberg procedures
applied to $A$ starting from an initial vector ${\bm q}_0$, then the Hessenberg matrix decomposition can
be uniformly written as
\begin{equation}
\begin{split}
AQ_m & = Q_mH_m + h_{m+1,m}{\bm q}_{m+1}{\bm e}^{\top}_m\\
     & = Q_{m+1}\bar{H}_m.
\end{split}
\label{eq2.3}
\end{equation}
Denote as $\sigma_m$ the smallest singular value of $H_{m+1,m} - [I_m;0]$, then ${\bm v}_m$ at Line~6 of Algorithm \ref{Alg2} is the corresponding right singular vector, and $Q_m{\bm v}_m$ is the approximate PageRank vector. The residual vector at each restarting cycle can be computed
as ${\bm r} = \sigma_mQ_{m + 1}{\bm u}_m$.
\label{them3.1}
\end{theorem}
\noindent{\bf Proof}. According to Eq. (\ref{eq2.3}) and Algorithm \ref{Alg2}, it follows
that
\begin{equation}
\begin{split}
{\bm r} & = A{\bm q}_m - {\bm q}_m = AQ_m{\bm v}_m - Q_m{\bm v}_m \\
        & = Q_{m+1}H_{m+1,m}{\bm v}_m - Q_m{\bm v}_m\\
        & = Q_{m+1}\left[H_{m+1,m} - \begin{pmatrix}
        I_m\\
        {\bm 0}
        \end{pmatrix}\right]{\bm v}_m\\
        & = \sigma_mQ_{m+1}{\bm u}_m,
\end{split}
\label{eq2.4}
\end{equation}
where ${\bm q}_m$ is an approximate PageRank vector. Thus, the assertion is verified. \hfill $\Box$


Below, we give the norms of the residual vectors computed by the Arnoldi and by the Hessenberg procedures, respectively:
\begin{equation}
\|{\bm r}\|_1 = \sigma_m\|Q_{m+1}{\bm u}_m\|_1,~~ {\rm where}~~
Q^{\top}_mQ_m:
\begin{cases}
= I_m\quad({\rm Arnoldi~process}),\\
\neq I_m\quad({\rm Hessenberg~process}),
\end{cases}
\label{eq2.5}
\end{equation}
and
\begin{equation}
\|{\bm r}\|_2 =
\begin{cases}
\sigma_m, & Q^{\top}_mQ_m = I_m \quad({\rm Arnoldi~process}),\\
\sigma_m\|Q_{m+1}{\bm u}_m\|_2,& Q^{\top}_mQ_m \neq I_m\quad({\rm Hessenberg~process}).
\end{cases}
\label{eq2.6}
\end{equation}

Although the 2-norm of the Arnoldi residual vectors is much cheaper than the 2-norm of the Hessenberg residual vectors, it should be noted that for PageRank computations it is generally recommended to use the 1-norm; refer, e.g.,
to \cite{Kamvar03,Wu2013x,Wu2010}. Therefore, the computational complexity of the stopping
criterion (at line~8 of Algorithm~\ref{Alg2}) is almost the same for both methods~\cite{Yin2010,Yin2012,Teramoto}.

Before we end this section, we provide estimates on the storage requirement and the computational complexity of the new algorithm, also compared against other popular methods.
\begin{table}[!htpb]
\centering
\caption{Memory requirement for running $k$ iterations of different PageRank algorithms.}
\begin{threeparttable}
\begin{tabular}{cccccccc}\toprule
Algorithm   & ${\rm dim}(n)$                                           & ${\rm dim}(k)$ & total \\
\midrule
\texttt{Power}       & ${\bm x}_k,~{\bm x}_{k-1}$                               & --             & $2n$  \\
\texttt{QE-Power}    & ${\bm x}_k,~{\bm x}_{k-1},~{\bm x}_{k-2},~{\bm x}_{k-3}$ & --             & $4n$  \\
\texttt{Inner-Outer} & ${\bm x},~{\bm y},~{\bm f}$                              & --             & $3n$  \\
\texttt{Arnoldi}     & $Q_k,~{\bm w}$                                           & $H_k$\tnote{1} & $(k+1)n + k^2/2 + 2k$   \\
\texttt{Hessenberg}  & $L_k,~{\bm u}$                                           & $H_k$          & $(k+1)n + k^2/2 + 2k$   \\
\texttt{A-Arnoldi}   & $Q_k,~{\bm w},~{\bm r}$                                  & $H_k$          & $(k + 2)n + k^2/2 + 2k$ \\
\bottomrule
\end{tabular}
\begin{tablenotes}
\item[1] To minimize the memory requirements, both $L_k$ and $\bar{H}_k$ can be written into the same array as $A$. Hence, the Hessenberg
process requires slightly less storage than the Arnoldi's procedure; refer to \cite{Sadok1999,Heyouni08} for discussions on this issue.
\end{tablenotes}
\end{threeparttable}
\label{tab1}
\end{table}
Table~\ref{tab1} shows the memory required in addition to $A$ for running $k$ iterations of the power method (referred to as~\texttt{Power} in the table), the power method with quadratic extrapolation (called as \texttt{QE-power}), the Arnoldi-type method (abbreviated as \texttt{Arnoldi}), the adaptively accelerated Arnoldi method (called as \texttt{A-Arnoldi}) and the Hessenberg-type method (abbreviated as \texttt{Hessenberg}). Here, ${\bm w}$, ${\bm x}$, ${\bm u}$ and ${\bm r}$ are intermediate working vectors used at the $k$-th step, and $Q_k$ denotes the $k$ orthonormal vectors in the modified Gram-Schmidt process. Analogously, $L_k$ denotes the $n\times k$ non-orthonormal matrix in the variant of the $LU$-like factorization process.

\begin{table}[t]
\centering
\caption{Computational cost of one cycle of different algorithms for computing PageRank.}
\begin{tabular}{cccccccc}
\hline
operation                                   & \texttt{Arnoldi}     & \texttt{Hessenberg}  & \texttt{A-Arnoldi}       \\
\hline
matrix-vector product (see Line 4 of Algorithm 1)                     & $2mN_z$     & $2mN_z$     & $2mN_z$         \\
inner product (e.g., see Line 6 of \texttt{Arnoldi})                             & $m(m + 1)n$ & 0           & $3m(m + 1)n/2$  \\
$({\bm u})_{{\bm p}(i)}$  (see Line 6 of Algorithm 1)                & 0           & $mn$         & 0               \\
\texttt{SAXPY}: ${\bm x} + \alpha{\bm y}$ (see Line 7 of Algorithm 1)          & $m(m + 1)n$ & $m(m + 1)n$ & $m(m+1)n$       \\
$\|{\bm u}\|_2$ or $\|{\bm u}\|_{\infty}$ (e.g., see Line 10 of Algorithm 1) & $2mn$       & $mn$        & $3mn$           \\
vector scaling (see Line 11 of Algorithm 1)                           & $mn$        & $mn$        & $mn$            \\

%
%

\hline
\end{tabular}
\label{tab2}
\end{table}

Table~\ref{tab2} shows the computational workloads required to execute
one cycle of each different
iterative algorithm. Here, $N_z$ represents the number of nonzero entries of matrix $A$.
In fact, both \texttt{Arnoldi} and \texttt{A-Arnoldi} for computing PageRank enjoy the similar
pseudo-code of Algorithm \ref{Alg2}, the only difference is to choose the Hessenberg, Arnoldi,
or generalized Arnoldi process at Line 4 of Algorithm~\ref{Alg2}. It implies that we need to
compare the cost of the Hessenberg process with both the Arnoldi and generalized Arnoldi procedures.
We can see that one cycle of the Hessenberg-type method is cheaper than for Arnoldi and
for the generalized Arnoldi methods. Thus its use can be computationally attractive for large PageRank computation.
The convergence performance of the
Hessenberg method is also superior to Arnoldi algorithms (i.e., \texttt{Arnoldi} and
\texttt{A-Arnoldi}), as proved numerically in the next section. Besides matrix-vector multiplications,
also the computation of vector norms and \texttt{SAXPY}\footnote{See the details from \url{https://developer.nvidia.com/blog/six-ways-saxpy/}.}~(which stands for ``Single-Precision
A$\cdot$X Plus Y" and is a combination of the scalar multiplication and vector addition) operations determines the total computational cost of these three algorithms. Overall, when $m$ increases, the cost of each cycle increases too but the total number of iterations  decreases. The optimal value of the restart parameter that minimizes the total solution time remains problem dependent, and this issue will be examined in our numerical experiments section.

\section{Numerical experiments}
\label{sec4}
In this section, numerical experiments are reported to illustrate the efficiency of
the Hessenberg-based PageRank algorithm presented in this paper also against other
popular PageRank algorithms, that are the conventional power method including its
variants with quadratic extrapolation~\cite{Kamvar03} and with linear extrapolation~\cite{Tan2017},
the Arnoldi-based PageRank method introduced in~\cite{Golub06}, and the adaptively accelerated Arnoldi method
\cite{Yin2012}\footnote{Numerical results with the  IDR($s$)-based PageRank method are omitted due to its unsatisfactory performance for large values of $s$ and $m$. However, the MATLAB code of the IDR(s)-based PageRank method is still included in our GitHub repository: \url{https://github.com/Hsien-Ming-Ku/PageRank-Hessenberg}
for testing purposes.}. The performance of these methods were assessed in terms of number of matrix-vector products (or, equivalently, number of iteration steps for the first three algorithms) and elapsed CPU time (in seconds) required to achieve convergence to a prescribed accuracy. Unless otherwise stated, the stopping criterion used in our runs was
\begin{equation*}
	\frac{\|A{\bm q} - {\bm q}\|_1}{\|{\bm q}\|_1} < tol ~= 10^{-8},
\end{equation*}
and all the algorithms were started from the initial vector ${\bm q}_0 = {\bm e}/\|{\bm e}\|_1$, where ${\bm e} = [1,\ldots, 1]^T$. According to Theorem \ref{them3.1}, the cost of implementing the above stopping criterion can be alleviated for both Arnoldi- and Hessenberg-type methods, since $A{\bm q} - {\bm q} = \sigma_m Q_{m + 1}{\bm u}_m$ and the computation of $\sigma_mQ_{m+1}{\bm u}_m$ is actually cheaper than that of $A{\bm q} - {\bm q}$, when $m$ is not large. In our experiments with the method denoted as \texttt{QE-power}, the quadratic extrapolation technique was applied every five iterations, following the observations made in~\cite{Kamvar03}. The experiments were run in MATLAB R2017b (64 bit) on a computer equipped with Intel Core i5-8250U processor (CPU 1.60$\sim$1.80 GHz), 8 GB of RAM using double precision floating point arithmetic with machine epsilon set equal to $10^{-16}$.

The matrix problems used in our runs are obtained from the \textit{SuiteSparse Matrix Collection}, which is available online
at \url{https://sparse.tamu.edu/}. In Table \ref{tab3}, we describe the characteristics of our test matrices,
including number of rows ($n$), number of nonzeros ($N_z$), number of zero columns ($zcol$), average nonzeros
of every row ($aN_z$) and density ($den$) which is defined as
\begin{equation*}
den = \frac{N_z}{n\times n}\times 100.
\end{equation*}
Here, the number of zero columns corresponds to the number of dangling nodes. The largest problem in our set has size $5,363,260$ and  $79,023,142$ nonzeros.

\begin{table}[!htpb]
\caption{The characteristic of test matrices.}
\vspace{1mm}
\centering
\begin{tabular}{ccccrrc}
\hline
Matrix ID                 & Matrix name           & Size      & Nonzeros   & $zcol$    & $aN_z$  & $den$                \\
\hline
\uppercase\expandafter{\romannumeral1} &\texttt{soc-Slashdot0902} & 82,168   & 948,464   & 3,727   & 11.543  & $1.405\times 10^{-2}$\\
\uppercase\expandafter{\romannumeral2} &\texttt{amazon0312}       & 400,727  & 3,200,440 & 12,353  & 7.987   & $1.993\times 10^{-3}$ \\
\uppercase\expandafter{\romannumeral3} & \texttt{amazon-2008}     & 735,323  & 5,158,388 & 88,557  & 7.015   & $9.540\times 10^{-4}$ \\
\uppercase\expandafter{\romannumeral4} & \texttt{wiki-Talk}       & 2,394,385 & 5,021,410 & 2,246,783 & 2.097    & $8.759\times 10^{-5}$ \\
\uppercase\expandafter{\romannumeral5} & \texttt{ljournal-2008}   & 5,363,260 & 79,023,142 & 545,626  & 14.734   & $2.747\times 10^{-4}$ \\
\hline
\end{tabular}
\label{tab3}
\end{table}

\subsection{Choice of the restart value $m$}
\begin{figure}[!htpb]
\centering
\includegraphics[width=0.48\textwidth]{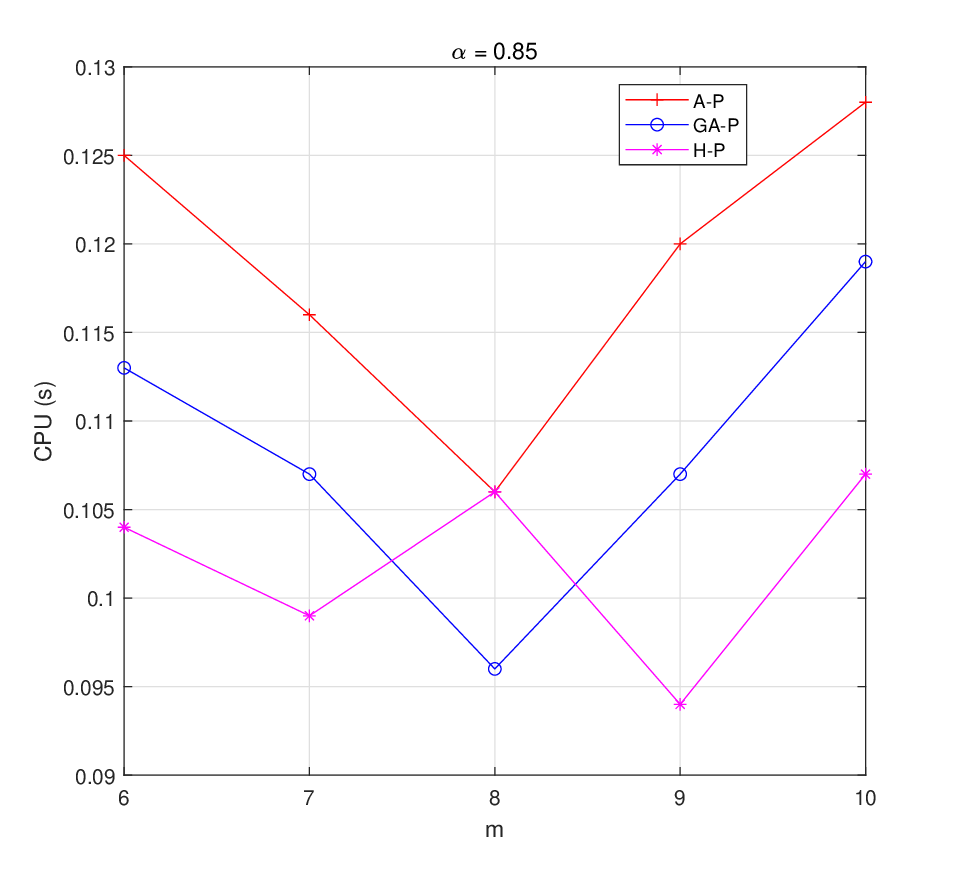}
\includegraphics[width=0.48\textwidth]{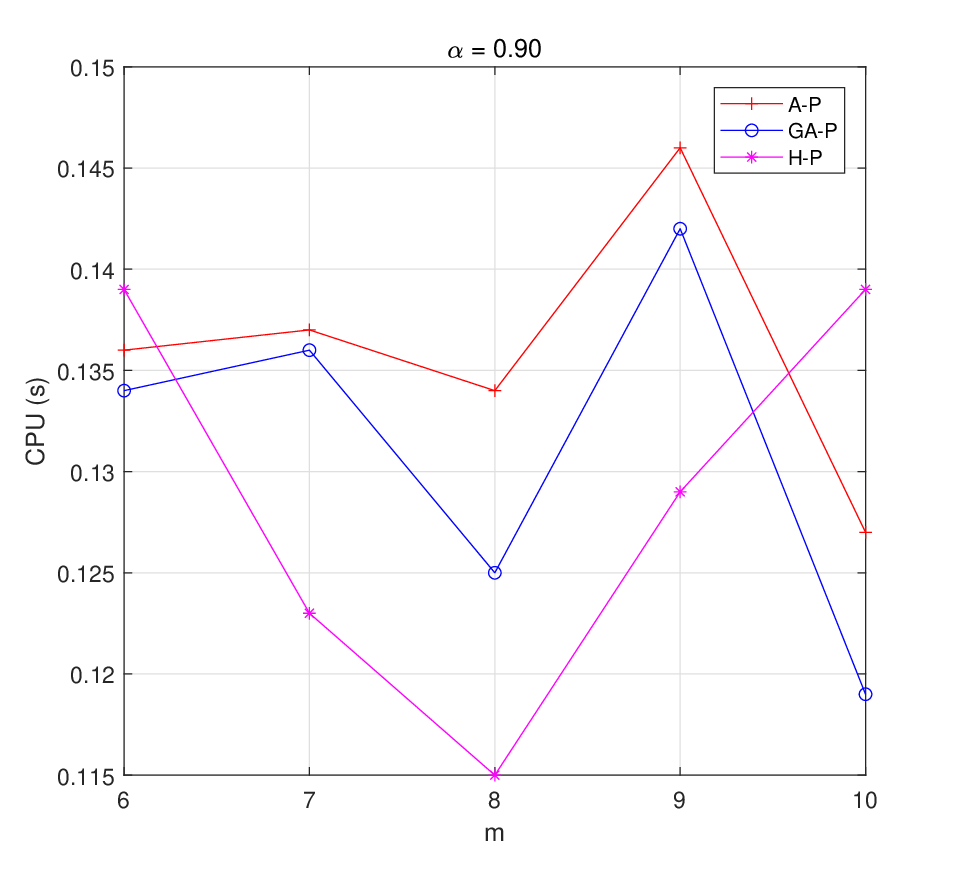}
\includegraphics[width=0.48\textwidth]{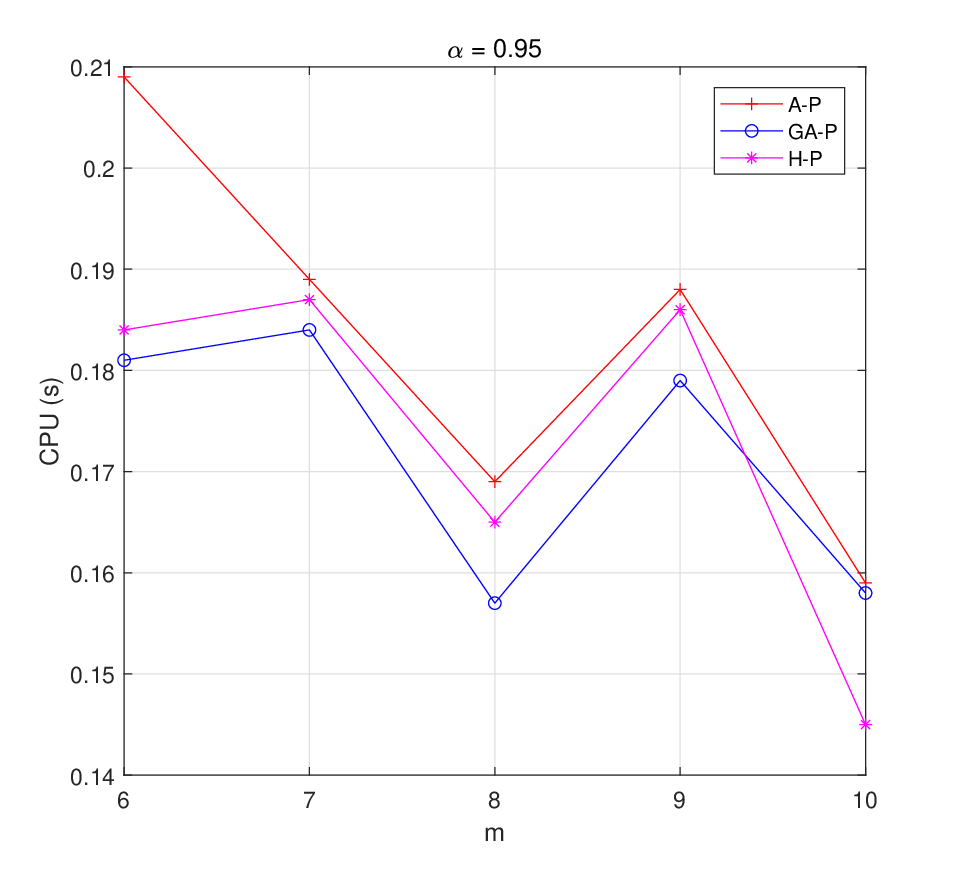}
\includegraphics[width=0.48\textwidth]{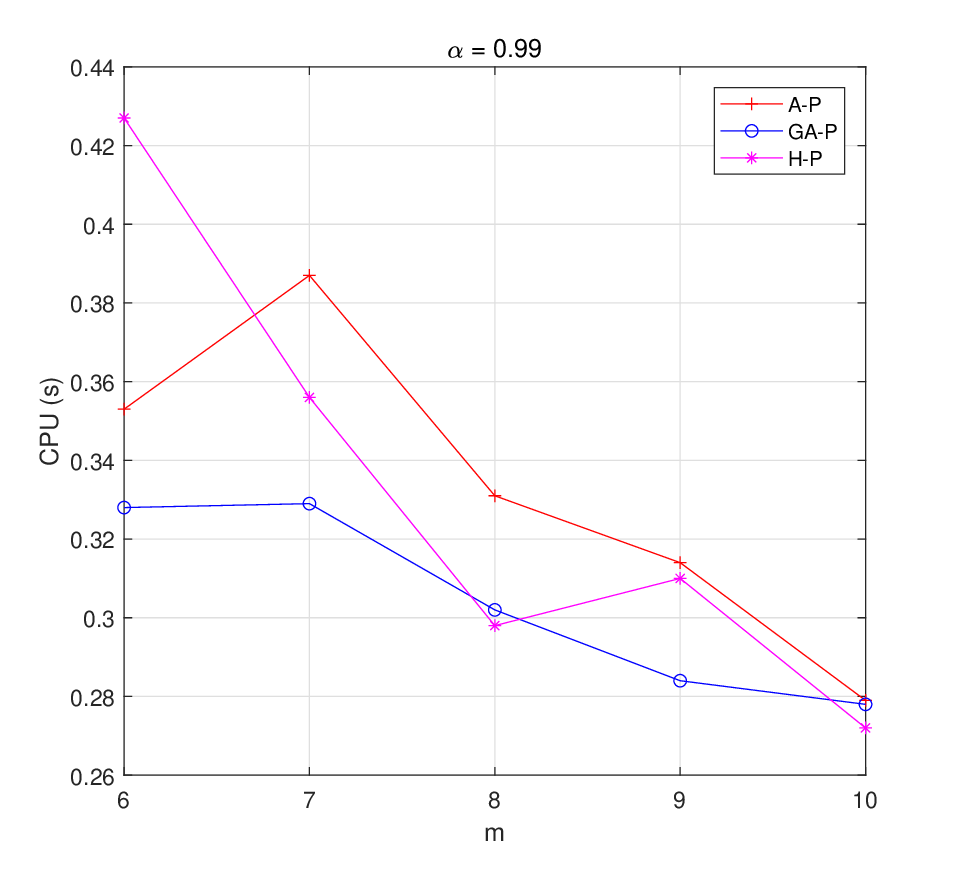}
\caption{Plot of the elapsed CPU time in seconds versus the restart number (i.e., $m$) for the
test problem `\texttt{soc-Slashdot0902}' using $tol = 10^{-7}$.}
\label{fig3}
\end{figure}

\begin{figure}[t]
\centering
\includegraphics[width=0.48\textwidth]{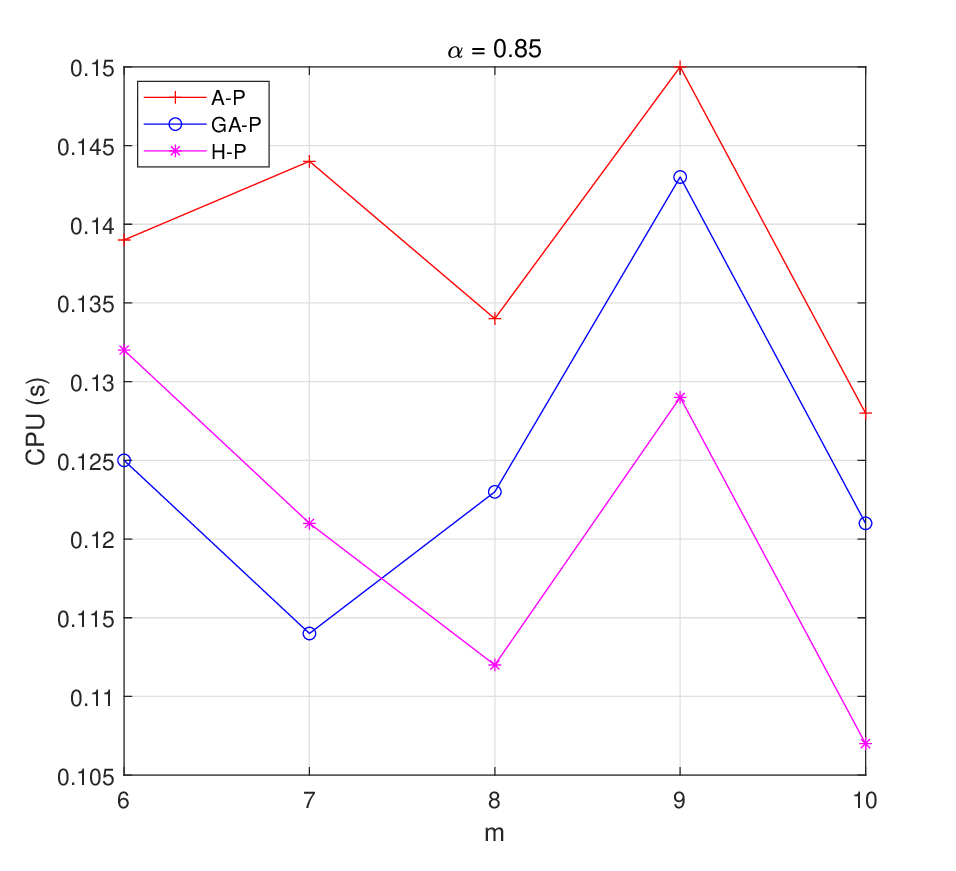}
\includegraphics[width=0.48\textwidth]{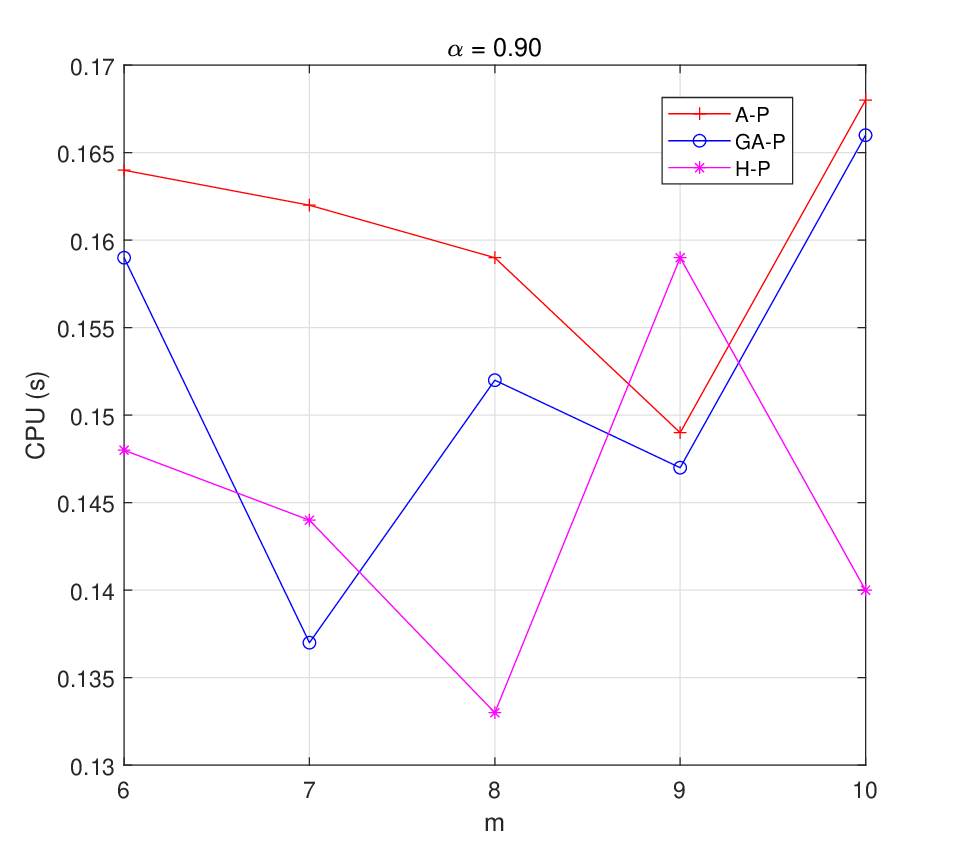}
\includegraphics[width=0.48\textwidth]{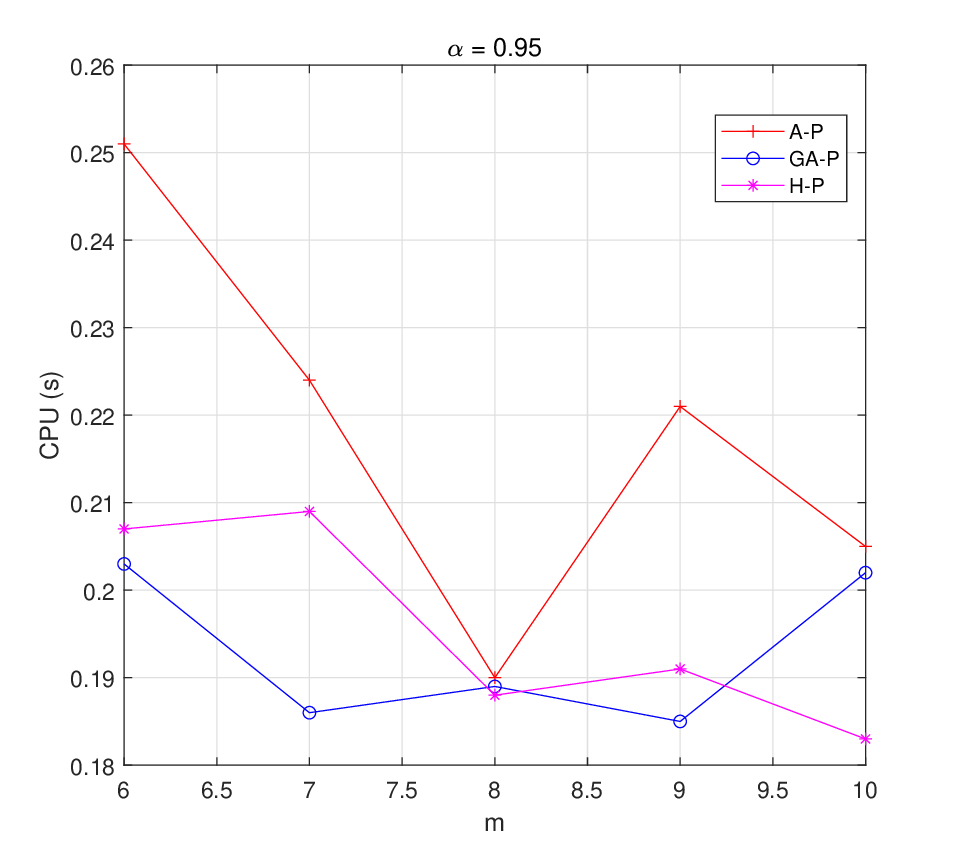}
\includegraphics[width=0.48\textwidth]{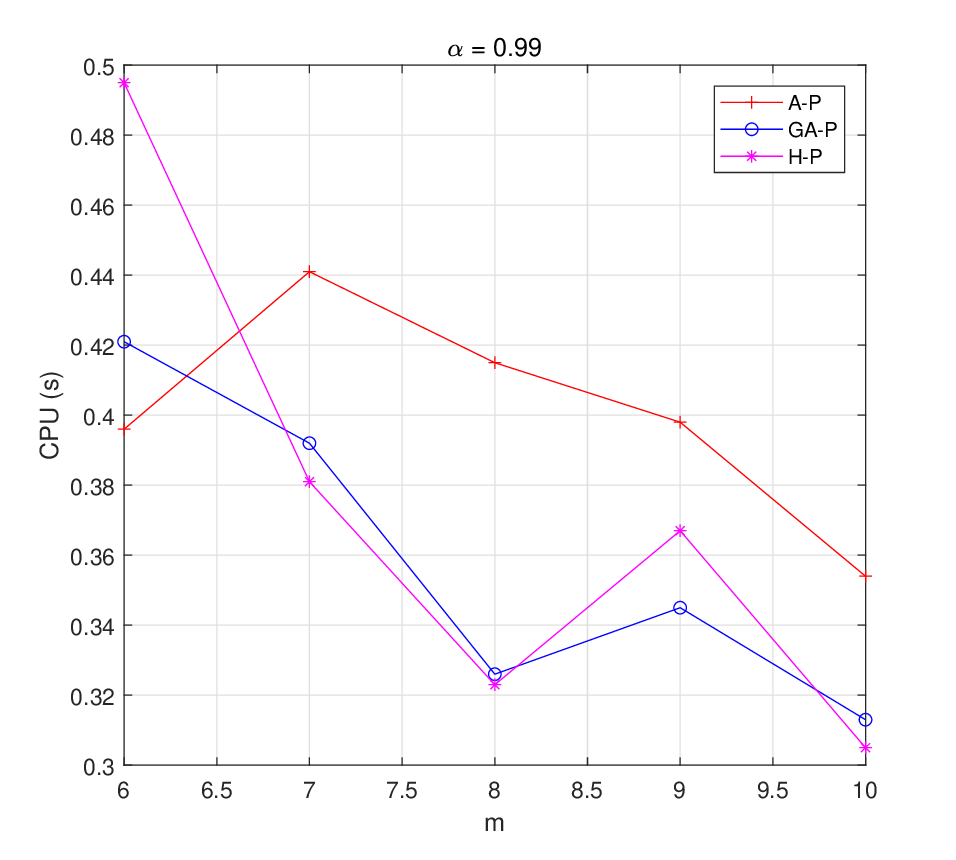}
\caption{Plot of the elapsed CPU time in seconds versus the restart number $m$ for the
test problem `\texttt{soc-Slashdot0902}' using $tol = 10^{-8}$.}
\label{fig4}
\end{figure}

First, we investigate the effect of the restart parameter $m$ on the convergence of the \texttt{Arnoldi (A-P)}, \texttt{A-Arnoldi (GA-P)} and \texttt{Hessenberg (H-P)} methods in terms of number of iterations and elapsed CPU time, since this parameter may noticeably affect the performance of Krylov subspace-based methods. The results are presented numerically in Tables \ref{tab4}--\ref{tab6}. In Figs.~\ref{fig3}--\ref{fig4}, for the test matrix `\texttt{soc-Slashdot0902}' we plot the curves showing the total CPU time versus $m$ for different damping factors and $tol$'s values.

\begin{table}[!htbp]\tabcolsep=5.5pt
\centering
\caption{Number of iterations required by the Arnoldi-, GArnoldi- and Hessenberg-based algorithms with
different restart numbers. (problem \texttt{soc-Slashdot0902} and $tol = 10^{-7}$).}
\vspace{1mm}
\begin{tabular}{lcccrcrcrcrcrccc}
\hline $\alpha$ &\multicolumn{3}{c}{$m=6$} &\multicolumn{3}{c}{$m=7$} &\multicolumn{3}
{c}{$m=8$}&\multicolumn{3}{c}{$m=9$}&\multicolumn{3}{c}{$m=10$}\\
[-2pt]\cmidrule(l{0.7em}r{0.7em}){2-4} \cmidrule(l{0.7em}r{0.6em}){5-7}\cmidrule(l{0.7em}r{0.7em}){8-10}
\cmidrule(l{0.7em}r{0.7em}){11-13}\cmidrule(l{0.7em}r{0.7em}){14-16}\\[-11pt]
& \texttt{A-P} &\texttt{GA-P} & $\texttt{H-P}$ &\texttt{A-P} &\texttt{GA-P} &$\texttt{H-P}$
&\texttt{A-P} &\texttt{GA-P} &$\texttt{H-P}$ &$\texttt{A-P}$ &$\texttt{GA-P}$&$\texttt{H-P}$ &$\texttt{A-P}$ &$\texttt{GA-P}$&$\texttt{H-P}$ \\
\hline
0.85 &5    &5   &5  &4    &4  &4    &3    &3   &4   &3    &3     &3  &3    &3     &3\\
0.90 &6    &6   &7  &5    &5  &5    &4    &4   &4   &4    &4     &4  &3    &3     &4\\
0.95 &9    &8   &9  &7    &7  &8    &5    &5   &6   &5    &5     &6  &4    &4     &4\\
0.99 &16   &15  &21 &15   &13 &15   &11   &10  &11  &9    &8     &10 &7    &7     &8\\
\hline
\end{tabular}
\label{tab4}
\end{table}
\begin{table}[!htbp]\tabcolsep=5.5pt
\centering
\caption{Number of iterations required by the Arnoldi-, GArnoldi- and Hessenberg-type algorithms with
different restart numbers (problem \texttt{soc-Slashdot0902} and $tol = 10^{-8}$).}
\vspace{1mm}
\begin{tabular}{lcccrcrcrcrcrccc}
\hline $\alpha$ &\multicolumn{3}{c}{$m=6$} &\multicolumn{3}{c}{$m=7$} &\multicolumn{3}{c}{$m=8$} &\multicolumn{3}
{c}{$m=9$}&\multicolumn{3}{c}{$m=10$}\\
[-2pt]\cmidrule(l{0.7em}r{0.7em}){2-4} \cmidrule(l{0.7em}r{0.6em}){5-7}\cmidrule(l{0.7em}r{0.7em}){8-10}
\cmidrule(l{0.7em}r{0.7em}){11-13}\cmidrule(l{0.7em}r{0.7em}){14-16}\\[-11pt]
& \texttt{A-P} &\texttt{GA-P} & $\texttt{H-P}$ &\texttt{A-P} &\texttt{GA-P} &$\texttt{H-P}$
&\texttt{A-P} &\texttt{GA-P} &$\texttt{H-P}$ &$\texttt{A-P}$ &$\texttt{GA-P}$&$\texttt{H-P}$&$\texttt{A-P}$ &$\texttt{GA-P}$&$\texttt{H-P}$ \\
\hline
0.85 &6   &5   &6   &5  &4   &5   &4    &4   &4    &4   &4   &4  &3    &3   &3     \\
0.90 &7   &7   &7   &6  &5   &6   &5    &5   &5    &4   &4   &5  &4    &4   &4    \\
0.95 &11  &9   &10  &8  &7   &9   &6    &6   &7    &6   &5   &6  &5    &5   &5    \\
0.99 &18  &19  &24  &17 &15  &16  &14   &11  &12   &11  &10  &12 &9    &8   &9     \\
\hline
\end{tabular}
\label{tab4yz}
\end{table}

According to the results reported in Tables \ref{tab4}--\ref{tab6}, the number
of iterations required to converge by these three algorithms tends to decrease
for higher restart numbers $m$, especially for larger damping factors. This
behaviour is expected because larger search spaces may provide better approximations.
On the other hand, the total solution time of the three methods is not significantly
reduced. As mentioned in Section~\ref{sec3}, and explained in Tables~\ref{tab1}--\ref{tab2},
the storage requirement and the computational cost of one Arnoldi and Hessenberg
cycles increase with $m$. However, it should be noted that \texttt{Hessenberg} is
more cost effective than both \texttt{Arnoldi} and \texttt{A-Arnoldi} for larger
$m$. In our numerical experiments, we choose the restart numbers {\color{red}equal
to $m = 8, 10$} due to memory constraints\footnote{In fact, the restart number is
often chosen in the PageRank literature as $m\leq 10$~\cite{Gu2017a,Teramoto,WuWanJin2012}.}, since the number of iterations and the total elapsed CPU time are still acceptable. It may be worth investigating techniques that can effectively reduce the dimension of the Krylov subspace for the \texttt{Hessenberg} method, e.g., by optimizing the choice of the starting vector~\cite{Astudillo} and by utilizing vector extrapolations \cite{Tan2017}, but this analyses is beyond the scope of this study.


\subsection{Effect of damping factors on the CPU time and the number of iterations}
For the five matrix problems listed in Table~\ref{tab3}, we report on the number of
matrix-vector products (\texttt{Mvp} in short) and the elapsed CPU time of the power
method, the power methods with  quadratic extrapolation and with linear extrapolation,
the Arnoldi-type method, the adaptively accelerated Arnoldi method and the Hessenberg-type
method for various values of the damping factor $\alpha$ ranging from $0.85$ to $0.99$.
\begin{sidewaystable}[!htpb]\small\tabcolsep=5.9pt
\begin{center}
\caption{Matrix-vector products and CPU time in seconds versus damping factors}
\begin{tabular}{llcccrcrcrcrcrcrcrcr}
\hline ID &$\alpha$ &\multicolumn{2}{c}{\texttt{Power}} &\multicolumn{2}{c}{\texttt{Power-Tan}} &\multicolumn{2}
{c}{\texttt{QE-Power}}&\multicolumn{4}{c}{\texttt{Arnoldi}}&\multicolumn{4}{c}{\texttt{A-Arnoldi}}&\multicolumn{4}{c}{\texttt{Hessenberg}}
\\
&   &   & && &&&\multicolumn{2}{c}{$m=8$}&\multicolumn{2}{c}{$m=10$}&\multicolumn{2}{c}{$m=8$}
&\multicolumn{2}{c}{$m=10$}&\multicolumn{2}{c}{$m=8$}&\multicolumn{2}{c}{$m=10$}
\\
[-2pt]\cmidrule(l{0.7em}r{0.7em}){3-4} \cmidrule(l{0.7em}r{0.6em}){5-6}\cmidrule(l{0.7em}r{0.7em}){7-8}
\cmidrule(l{0.7em}r{0.7em}){9-10}\cmidrule(l{0.7em}r{0.7em}){11-12}\cmidrule(l{0.7em}r{0.7em}){13-14}
\cmidrule(l{0.7em}r{0.7em}){15-16}\cmidrule(l{0.7em}r{0.7em}){17-18}\cmidrule(l{0.7em}r{0.7em}){19-20}\\[-11pt]
& &\texttt{Mvp} &\texttt{CPU} & $\texttt{Mvp}$ &\texttt{CPU} &\texttt{Mvp} &
$\texttt{CPU}$&\texttt{Mvp} &\texttt{CPU} &$\texttt{Mvp}$ &$\texttt{CPU}$ &$\texttt{Mvp}$ &$\texttt{CPU}$
&$\texttt{Mvp}$ &$\texttt{CPU}$&$\texttt{Mvp}$ &$\texttt{CPU}$&$\texttt{Mvp}$ &$\texttt{CPU}$\\
\hline
\uppercase\expandafter{\romannumeral1}
   &0.85 &82   &0.216   &83   &0.228   &52   &0.168   &32   &0.152  &30   &0.136   &32     &0.122
         &30   &0.132   &32   &0.122   &30   &{\bf 0.118}  \\
   &0.90 &125  &0.324   &126  &0.333   &52   &0.159   &40   &0.186  &40   &0.187   &40     &0.159
         &40   &0.163   &40   &{\bf 0.144}   &40   &0.146  \\
   &0.95 &241  &0.628   &243  &0.634   &98   &0.298   &48   &0.192  &50   &0.206   &48     &0.195
         &50   &0.202   &56   &0.187   &50   &{\bf 0.176}  \\
   &0.99 &824  &2.089   &832  &2.105   &169  &0.504   &112  &0.419  &90   &0.350   &88     &0.327
         &80   &0.308   &96   &0.319   &90   &{\bf 0.303}  \\
\hline
\uppercase\expandafter{\romannumeral2}
   &0.85 &80   &1.032   &80   &1.035   &66   &1.014   &40   &0.951  &40   &0.943   &40     &0.996
         &40   &1.089   &48   &0.854   &40   &{\bf 0.726}   \\
   &0.90 &120  &1.553   &121  &1.562   &98   &1.522   &56   &1.253  &40   &1.187   &48     &1.166
         &40   &1.292   &56   &0.993   &40   &{\bf 0.918}   \\
   &0.95 &234  &3.033   &235  &3.041   &152  &2.422   &80   &1.779  &56   &1.681   &72     &1.732
         &70   &1.802   &88   &1.554   &80   &{\bf 1.393}   \\
   &0.99 &1071 &13.80   &1076 &14.16   &391  &6.196   &224  &4.917  &190  &4.404   &176    &4.212
         &160  &4.095   &216  &{\bf 3.798}   &250  &4.387   \\
\hline
\uppercase\expandafter{\romannumeral3}
   &0.85 &79   &1.782   &80   &1.821   &64   &1.779   &40   &1.633  &40   &1.704   &40     &1.761
         &40   &1.896   &40   &{\bf 1.315}   &40   &1.321   \\
   &0.90 &119  &2.782   &120  &2.839   &95   &2.659   &48   &1.981  &50   &2.304   &48     &2.130
         &50   &2.376   &56   &1.828   &50   &{\bf 1.627}   \\
   &0.95 &235  &5.248   &236  &5.289   &152  &4.215   &72   &2.889  &70   &3.034   &72     &3.152
         &70   &3.283   &72   &{\bf 2.345}   &80   &2.632   \\
   &0.99 &1082 &24.85   &1084 &25.02   &411  &11.79   &208  &8.428  &180  &7.779   &152    &6.817
         &150  &7.107   &176  &{\bf 5.843}   &190  &6.358   \\
\hline
\uppercase\expandafter{\romannumeral4}
   &0.85 &67   &4.798   &67    &4.819  &52   &4.727   &32   &4.321  &30   &4.403   &32     &4.683
         &30   &4.808   &32    &{\bf 3.327}  &40   &4.264   \\
   &0.90 &104  &7.344   &104   &7.464  &52   &4.770   &40   &5.352  &40   &5.838   &40     &5.948
         &40   &6.266   &40    &{\bf 4.217}  &40   &4.299   \\
   &0.95 &214  &15.11   &214   &15.27  &84   &7.327   &56   &7.413  &60   &8.736   &48     &7.108
         &60   &9.532   &64    &6.768  &50   &{\bf 5.335}   \\
   &0.99 &1066 &20.21   &1068  &75.96  &175  &{\bf 15.57}   &156  &20.21  &150  &21.47   &112    &16.24
         &110  &17.08   &160   &16.69  &150  &15.93   \\
\hline
\uppercase\expandafter{\romannumeral5}
   &0.85 &80   &36.59   &79    &36.47  &61   &30.61   &48   &28.53  &50   &30.88   &48     &29.79
         &50   &32.38   &56    &29.44  &50   &{\bf26.87}   \\
   &0.90 &122  &54.77   &122   &55.12  &92   &45.32   &64   &38.26  &60   &36.66   &64     &39.21
         &60   &38.65   &72    &37.62  &60   &{\bf31.82}   \\
   &0.95 &247  &110.9   &244   &109.5  &148  &71.61   &104  &61.85  &90   &55.92   &96     &59.05
         &90   &57.84   &104   &54.59  &90   &{\bf47.97}   \\
   &0.99 &1153 &794.7   &1061  &529.1  &630  &322.6   &352  &214.7  &240  &145.8   &216    &132.3
         &190  &121.4   &248   &127.7  &220  &{\bf115.9} \\
\hline
\end{tabular}
\label{tab6}
\end{center}
\end{sidewaystable}

We can see from the results of Table~\ref{tab6} that the power method accelerated by quadratic extrapolation outperforms the conventional power method and its linearly extrapolated variant, while in most cases our Hessenberg-based solver is the fastest method in terms of elapsed CPU time, with the only exception for matrix
`\uppercase\expandafter{\romannumeral4}' using $\alpha = 0.99$. It can be observed that \texttt{Arnoldi} is more cost-effective than \texttt{A-Arnoldi} at equal number of \texttt{Mvp}, especially for large problems. This behaviour is in agreement with the cost analysis presented in Table~\ref{tab2}. Except these few cases, the \texttt{A-Arnoldi} method is still attractive to consider. On the other hand, one observes from
Table~\ref{tab6} that the numerical behaviors of the \texttt{Arnoldi}, \texttt{A-Arnoldi}
and \texttt{Hessenberg} algorithms relies on the choice of $m$ and $\alpha$.
For example, when $m$ is small, say $m=8$, these three algorithms are only slightly better than the \texttt{Power},
\texttt{Power-Tan} and \texttt{QE-Power} methods. However, as $m$ and $\alpha$ increase, their
improvements become gradually more significant. Unlike the \texttt{Arnoldi}, \texttt{A-Arnoldi},
and \texttt{Hessenberg}, the \texttt{Power}, \texttt{Power-Tan} and \texttt{QE-Power} methods
are simple and their main computational cost is the evaluation of matrix-vector products.
These characteristics often make them still feasible for computing PageRank when the damping factor $\alpha$ is not large.

In addition, it is interesting to mention that the Hessenberg-type method often needs more \texttt{Mvp}'s for convergence than
both \texttt{Arnoldi} and \texttt{A-Arnoldi}, whereas the total CPU time of \texttt{Hessenberg}
is still less. This is because the \texttt{Hessenberg} uses the cheap similarity transformations to reduce
the large matrix into the Hessenberg form, whereas the latter two methods use expensive (weighted) unitary transformations.

\FloatBarrier

\section{Conclusions}
\label{sec5}
In this paper, we proposed a novel approach for computing the PageRank problem.
The proposed method has lower computational cost than both \texttt{Arnoldi} and \texttt{A-Arnoldi} to find the approximate PageRank vector, thus it can afford to use
higher dimensional
Krylov subspaces. Extensive numerical experiments are
reported to illustrate the efficiency of the proposed method also compared to other state-of-the-art matrix solvers for this problem class, especially when the damping factor is large. Hence, we conclude that the $\texttt{Hessenberg}$ method as well as the \texttt{Arnoldi} and \texttt{A-Arnoldi} methods can be useful computational tools for practical large-scale PageRank computations.

Future research will focus on the theory of the Hessenberg process and
the convergence of the Hessenberg-type algorithm is still required to
be further analyzed. In addition, it is interesting to study how to reduce
the restart number $m$ and improve the convergence speed of our methods.
Moreover, the proposed method can be extended to compute the more general
Markov chains \cite{Freund94,Gleich15}, e.g., in ProteinRank and CiteRank.
\section*{Acknowledgements}
{\em The authors would like to thank Prof. Zhongxiao Jia for his comments about the strategy used in the refined
Arnoldi algorithm. Meanwhile, the authors are grateful to Dr. Reinaldo Astudillo (ASML Holding N.V.)
for his kind suggestions about executing the IDR-based Hessenberg decompositions used in Section \ref{sec2.2}.
This research is supported by NSFC (11601323 and 11801463), the Applied Basic Research Program
of Sichuan Province (2020YJ0007) and the research grants MYRG2018-00025-FST, MYRG2020-00208-FST from University of Macau. The last author is member of the Gruppo Nazionale per il Calcolo Scientifico (GNCS) of the Istituto Nazionale di Alta Matematica (INdAM) and his work was partially supported by INdAM-GNCS under Progetti di Ricerca 2020.}

\def\refname

{\large \bfseries References}
\begin{thebibliography}{1}
\bibitem{Page1998}
L. Page, S. Brin, R. Motwani, T. Winograd, The PageRank citation ranking:
Bringing order to the web, {\em Technical Report No. 1999-66}, Stanford
InfoLab., Jan. 29, 1999, 17 pages. Available online at: \url{http://ilpubs.stanford.edu:8090/422/}.

\bibitem{Kamvar03}
S.D. Kamvar, T.H. Haveliwala, C.D. Manning, G.H. Golub, Extrapolation methods
for accelerating PageRank computations, in: {\em WWW '03 Proceedings of the
12th international conference on World Wide Web}, Budapest, Hungary, May 20-24,
2003, ACM New York, NY (2003): 261-270. DOI: \href{https://doi.org/10.1145/775152.775190}{10.1145/775152.775190}.

\bibitem{Kamvar04x}
S. Kamvar, T. Haveliwala, G. Golub, Adaptive methods for the computation of
PageRank, {\em Linear Algebra Appl.}, 386 (2004): 51-65.

\bibitem{Langv05}
A.N. Langville, C.D. Meyer, Deeper inside PageRank, {\em Internet Math.}, 1(3) (2005): 335-380.

\bibitem{Langville05}
A.N. Langville, C.D. Meyer, A survey of eigenvector methods of web information
retrieval, {\em SIAM Rev.}, 47(1) (2005): 135-161.

\bibitem{Berkhin05}
P. Berkhin, A survey on PageRank computing, {\em Internet Math.}, 2(1)
(2005): 73-120.

\bibitem{LanMey2006}
A.N. Langville, C.D. Meyer, {\em Google's PageRank and Beyond: The Science of
Search Engine Rankings}, Princeton University Press, Princeton, NJ (2006).

\bibitem{Gleich15}
D.F. Gleich, PageRank beyond the web, {\em SIAM Rev.}, 57(3) (2015): 321-363.

\bibitem{Bryan06}
K. Bryan, T. Leise, The $25,000,000,000$ eigenvector: The linear algebra behind Google,
{\em SIAM Rev.}, 48(3) (2006): 569-581.

\bibitem{Cicone10}
A. Cicone, S. Serra-Capizzano, Google PageRanking problem: The model and the
analysis, {\em J. Comput. Appl. Math.}, 234(11) (2010): 3140-3169.

\bibitem{Avrach07}
K. Avrachenkov, N. Litvak, D. Nemirovsky, N. Osipova, Monte Carlo methods in
PageRank computation: when one iteration is sufficient, {\em SIAM J. Numer.
Anal.}, 45(2) (2007): 890-904.

\bibitem{Liu2015}
W. Liu, G. Li, J. Cheng, Fast PageRank approximation by adaptive sampling, {\em
Knowl. Inf. Syst.}, 42(1) (2015): 127-146.

\bibitem{Tan2017}
X. Tan, A new extrapolation method for PageRank computations, {\em J. Comput.
Appl. Math.}, 313 (2017): 383-392.

\bibitem{Langville06}
A.N. Langville, C.D. Meyer, A reordering for the PageRank problem, {\em SIAM
J. Sci. Comput.}, 27(6) (2006): 2112-2120.

\bibitem{Lin2009}
Y. Lin, X. Shi, Y. Wei, On computing PageRank via lumping the Google matrix,
{\em J. Comput. Appl. Math.}, 224(2) (2009): 702-708.

\bibitem{Gleich10}
D.F. Gleich, A.P. Gray, C. Greif, T. Lau, An inner-outer iteration for
computing PageRank, {\em SIAM J. Sci. Comput.}, 32(1) (2010): 349-371.

\bibitem{Heyouni98}
M. Heyouni, H. Sadok, On a variable smoothing procedure for Krylov subspace methods, {\em Linear Algebra Appl.},
268 (1998): 131-149.

\bibitem{Saad2010}
Y. Saad, {\em Numerical Methods for Large Eigenvalue Problems} (Revised Ed.),
SIAM, Philadelphia, PA (2011).

\bibitem{Jia1997}
Z. Jia, Refined iterative algorithms based on Arnoldi's process
for large unsymmetric eigenproblems, {\em Linear Algebra Appl.},
259 (1997): 1-23.

\bibitem{Jia1999}
Z. Jia, Polynomial characterizations of the approximate eigenvectors
by the refined Arnoldi method and an implicitly restarted refined
Arnoldi algorithm, {\em Linear Algebra Appl.}, 287(1-3) (1999):
191-214.

\bibitem{Jia2000}
Z. Jia, A refined subspace iteration algorithm for large sparse
eigenproblems, {\em Appl. Numer. Math.}, 32(1) (2000): 35-52.

\bibitem{Golub06}
G.H. Golub, C. Greif, An Arnoldi-type algorithm for computing page
rank, {\em BIT}, 46(4) (2006): 759-771.

\bibitem{Wu2007}
G. Wu, Y. Wei, A Power-Arnoldi algorithm for computing PageRank, {\em Numer.
Linear Algebra Appl.}, 14(7) (2007): 521-546.

\bibitem{Yin2010}
G.-J. Yin, J.-F. Yin, On Arnoldi method accelerating PageRank cmputations,
in: {\em Web Information Systems and Mining. WISM 2010 (F.-L. Wang, Z. Gong,
X. Luo, J. Lei, eds.)}, Lecture Notes in Computer Science, vol 6318, Springer,
Berlin, Heidelberg (2010): 378-385. DOI: \href{https://doi.org/10.1007/978-3-642-16515-3_47}{10.1007/978-3-642-16515-3\_47}.

\bibitem{Wu2013x}
G. Wu, Y. Zhang, Y. Wei, Accelerating the Arnoldi-type algorithm for the
PageRank problem and the ProteinRank problem, {\em J. Sci. Comput.}, 57(1)
(2013): 74-104.

\bibitem{Wu2010}
G. Wu, Y. Wei, An Arnoldi-extrapolation algorithm for computing PageRank,
{\em J. Comput. Appl. Math.}, 234(11) (2010): 3196-3212.

\bibitem{Gu2017a}
C. Gu, W. Wang, An Arnoldi-Inout algorithm for computing PageRank problems, {\em
J. Comput. Appl. Math.}, 309 (2017): 219-229.

\bibitem{Yin2012}
J.-F. Yin, G.-J. Yin, M. Ng, On adaptively accelerated Arnoldi method for
computing PageRank, {\em Numer. Linear Algebra Appl.}, 19(1) (2012):
73-85.

\bibitem{Freund94}
R.W. Freund, M. Hochbruck, On the use of two QMR algorithms for solving singular
systems and applications in Markov chain modeling, {\em Numer. Linear Algebra Appl.},
1(4) (1994): 403-420.

\bibitem{Teramoto}
K. Teramoto, T. Nodera, A note on Lanczos algorithm for computing PageRank, in:
{\em Forging Connections between Computational Mathematics and Computational
Geometry (K. Chen, A. Ravindran, eds.)}, Springer Proceedings in Mathematics
\& Statistics, Vol. 124, Springer, Cham, Switzerland (2016): 25-33. DOI:
\href{https://doi.org/10.5176/2251-1911_CMCGS14.15_3}{10.5176/2251-1911\_CMCGS14.15\_3}.

\bibitem{WuWanJin2012}
G. Wu, Y.-C. Wang, X.-Q. Jin, A preconditioned and shifted GMRES algorithm for the PageRank problem with multiple
damping factors, \emph{SIAM J. Sci. Comput.}, 34(5) (2012): A2558-A2575.

\bibitem{Wu2010y}
G. Wu, Y. Wei, Arnoldi versus GMRES for computing pageRank: A theoretical
contribution to google's pageRank problem, {\em ACM Trans Inf. Syst.}, 28(3)
(2010): 11. DOI: \href{https://doi.org/10.1145/1777432.1777434}{10.1145/1777432.1777434}.

\bibitem{Hessen40}
K. Hessenberg, {\em BehandLung Linearer Eigenwertaufgaben Mit Hilfe Der
Hamilton-Cayleyschen Gleichung}, Numerische Verfahren, Bericht 1, Institut
f\"{u}r Praktische Mathematik (IPM), Technische Hochschule Darmstadt (1940).
The scanned report and a biographical sketch of Karl Hessenberg's life are
available at \url{http://www.hessenberg.de/karl1.html}.

\bibitem{Wilkinson}
J.H. Wilkinson, {\em The Algebraic Eigenvalue Problem}, Clarendon Press,
Oxford, UK (1965).

\bibitem{Sadok1999}
H. Sadok, CMRH: A new method for solving nonsymmetric linear systems
based on the Hessenberg reduction algorithm, {\em Numer. Algorithms},
20(4) (1999), pp. 303-321.

\bibitem{Stephens99}
D. Stephens, {\em ELMRES: An Oblique Projection Method To Solve Sparse
Non-Symmetric Linear Systems (Ph.D dissertation)}, Florida Institute of
Technology, Melbourne, USA (1999). Available online at
\url{http://ncsu.edu/hpc/Documents/Publications/gary_howell/stephens.pdf}.

\bibitem{House2010}
A.S. Householder, F.L. Bauer, On certain methods for expanding the characteristic
polynomial, {\em Numer. Math.}, 1(1) (1959): 29-37.

\bibitem{Sadok2012}
H. Sadok, D.B. Szyld, A new look at CMRH and its relation to GMRES,
{\em BIT}, 52(2) (2012): 485-501.

\bibitem{Heyouni08}
M. Heyouni, H. Sadok, A new implementation of the CMRH method for solving
dense linear systems, {\em J. Comput. Appl. Math.}, 213(2) (2008): 387-399.

\bibitem{Zhang14x}
K. Zhang, C. Gu, Flexible global generalized Hessenberg methods for linear systems
with multiple right-hand sides, {\em J. Comput. Appl. Math.}, 263 (2014): 312-325.

\bibitem{Heyouni01}
M. Heyouni, The global Hessenberg and CMRH methods for linear systems with
multiple right-hand sides, {\em Numer. Algorithms}, 26(4) (2001): 317-332.

\bibitem{Heyouni05x}
M. Heyouni, A. Essai, Matrix Krylov subspace methods for linear systems
with multiple right-hand sides, {\em Numer. Algorithms}, 40(2) (2005): 137-156.

\bibitem{Amini18x}
S. Amini, F. Toutounian, M. Gachpazan, The block CMRH method for solving
nonsymmetric linear systems with multiple right-hand sides, {\em J. Comput.
Appl. Math.}, 337 (2018): 166-174.

\bibitem{Amini18y}
S. Amini, F. Toutounian, Weighted and flexible versions of block CMRH method
for solving nonsymmetric linear systems with multiple right-hand sides, {\em
Comput. Math. Appl.}, 76(8) (2018): 2011-2021.

\bibitem{Gu2018JCAM}
X.-M. Gu, T.-Z. Huang, G. Yin, B. Carpentieri, C. Wen, L. Du, Restarted
Hessenberg method for solving shifted nonsymmetric linear systems, {\em
J. Comput. Appl. Math.}, 331 (2018): 166-177.

\bibitem{Gu2018X}
X.-M. Gu, T.-Z. Huang, B. Carpentieri, A. Imakura, K. Zhang, L. Du, Efficient
variants of the CMRH method for solving a sequence of multi-shifted non-Hermitian
linear systems simultaneously, {\em J. Comput. Appl. Math.}, 375 (2020): 112788. DOI: \href{https://doi.org/10.1016/j.cam.2020.112788}{10.1016/j.cam.2020.112788}.

\bibitem{Ramezani18}
Z. Ramezani, F. Toutounian, Extended and rational Hessenberg methods for the
evaluation of matrix functions, {\em BIT}, 59(2) (2019): 523-545.




\bibitem{Addam2017}
M. Addam, M. Heyouni, H. Sadok, The block Hessenberg process for matrix equations,
{\em Electron. Trans. Numer. Anal.}, 46 (2017): 460-473.

\bibitem{Heyouni19}
M. Heyouni, F. Saberi-Movahed, A. Tajaddini, On global Hessenberg based methods
for solving Sylvester matrix equations, {\em Comput. Math. Appl.}, 77(1) (2019): 77-92.

\bibitem{Businger}
P.A. Businger, Reducing a matrix to Hessenberg form, {\em Math. Comp.},
23(108) (1969): 819-821.

\bibitem{Heyouni99x}
M. Heyouni, Newton Generalized Hessenberg method for solving nonlinear
systems of equations, {\em Numer. Algorithms}, 21(1-4) (1999): 225-246.

\bibitem{Astudillo}
R. Astudillo, M.B. van Gijzen, A restarted Induced Dimension Reduction
method to approximate eigenpairs of large unsymmetric matrices, {\em J.
Comput. Appl. Math.}, 296 (2016): 24-35.
%

\bibitem{Gutknecht}
M.H. Gutknecht, J.-P.M. Zemke, Eigenvalue computations based on IDR,
{\em SIAM J. Matrix Anal. Appl.}, 34(2) (2013): 283-311.
%
\end{thebibliography}
\end{document}